\documentclass[11pt]{article}
\usepackage{amsmath,amsfonts}
\usepackage{verbatim}
\usepackage{latexsym}
\usepackage{graphicx}
\usepackage{float}
\usepackage{color}
\usepackage{cite}
\usepackage[hidelinks]{hyperref}
\usepackage{amssymb}
\usepackage{a4}
\allowdisplaybreaks
\makeatletter
\@addtoreset{equation}{section}

\begin{document}

\newcommand{\E}{\mathbb{E}}
\newcommand{\PP}{\mathbb{P}}
\newcommand{\RR}{\mathbb{R}}

\newtheorem{theorem}{Theorem}[section]
\newtheorem{lemma}{Lemma}[section]
\newtheorem{coro}{Corollary}[section]
\newtheorem{defn}{Definition}[section]
\newtheorem{assp}{Assumption}[section]
\newtheorem{expl}{Example}[section]
\newtheorem{prop}{Proposition}[section]
\newtheorem{rmk}{Remark}[section]

\newcommand\tq{{\scriptstyle{3\over 4 }\scriptstyle}}
\newcommand\qua{{\scriptstyle{1\over 4 }\scriptstyle}}
\newcommand\hf{{\textstyle{1\over 2 }\displaystyle}}
\newcommand\hhf{{\scriptstyle{1\over 2 }\scriptstyle}}

\newcommand{\proof}{\noindent {\it Proof}. }
\newcommand{\eproof}{\hfill $\Box$} 

\def\a{\alpha} \def\g{\gamma}
\def\e{\varepsilon} \def\z{\zeta} \def\y{\eta} \def\o{\theta}
\def\vo{\vartheta} \def\k{\kappa} \def\l{\lambda} \def\m{\mu} \def\n{\nu}
\def\x{\xi}  \def\r{\rho} \def\s{\sigma}
\def\p{\phi} \def\f{\varphi}   \def\w{\omega}
\def\q{\surd} \def\i{\bot} \def\h{\forall} \def\j{\emptyset}

\def\be{\beta} \def\de{\delta} \def\up{\upsilon} \def\eq{\equiv}
\def\ve{\vee} \def\we{\wedge}

\def\F{{\cal F}}
\def\T{\tau} \def\G{\Gamma}  \def\D{\Delta} \def\O{\Theta} \def\L{\Lambda}
\def\X{\Xi} \def\S{\Sigma} \def\W{\Omega}
\def\M{\partial} \def\N{\nabla} \def\Ex{\exists} \def\K{\times}
\def\V{\bigvee} \def\U{\bigwedge}

\def\1{\oslash} \def\2{\oplus} \def\3{\otimes} \def\4{\ominus}
\def\5{\circ} \def\6{\odot} \def\7{\backslash} \def\8{\infty}
\def\9{\bigcap} \def\0{\bigcup} \def\+{\pm} \def\-{\mp}
\def\la{\langle} \def\ra{\rangle}

\def\tl{\tilde}
\def\trace{\hbox{\rm trace}}
\def\diag{\hbox{\rm diag}}
\def\for{\quad\hbox{for }}
\def\refer{\hangindent=0.3in\hangafter=1}

\newcommand\wD{\widehat{\D}}
\title{
\bf The adaptive EM schemes for McKean-Vlasov SDEs with common noise in finite and infinite horizons
}

\author{
{\bf ${Hu\ Liu}^{1}$, ${Shuaibin\ Gao}^{2}$, ${Junhao\ Hu}^{2}$\thanks{The corresponding author. Email: junhaohu74@163.com}},
\\1 School of Mathematics and Information Science, North Minzu University,\\
Yinchuan, 750021,  China.\\
2 School of Mathematics and Statistics, 
South-Central Minzu University, \\
Wuhan, 430074, China.\\
}

\date{}

\maketitle

\begin{abstract}
\indent This paper is dedicated to investigating the adaptive Euler-Maruyama (EM) schemes for the approximation of McKean-Vlasov stochastic differential equations (SDEs) with common noise. When the drift and diffusion coefficients both satisfy the superlinear growth conditions, the $L^p$ convergence rates in finite and infinite horizons are revealed, which reacts to the particle number and step size. Subsequently, there is an illustration of the theory results by means of two numerical examples.

\medskip \noindent
{\small\bf Key words.}  McKean-Vlasov stochastic differential equations; particle systems; adaptive EM scheme; common noise; convergence
 \par \noindent
\end{abstract}

\section{Introduction}\indent
\indent  The McKean-Vlasov SDE (introduced by McKean in \cite{M66}) for a $d$-dimensional valued process  $X=(X_t)_{t\geq0}$, is a type of SDEs whose coefficients depend on both state variables and the distributions of state, that is,
\begin{equation}\label{eq1.1} 
	dX_t =b(X_t,{\cal L}_{X_t})dt + \sigma(X_t,{\cal L}_{X_t}) dW_t, 
\end{equation}
 where $W=(W_t)_{t\geq0}$ represents an $m$-dimensional standard Brownian motion and ${\cal L}_{X_t}$ denotes the marginal law of the process $X$ at time $t$. The initial condition is given by the $\RR^d$-valued random variable $X_0=\xi\in{L}_0^{\tilde{p}}(\mathbb{R}^d),$ where ${L}_0^{\tilde{p}}(\mathbb{R}^d)$ is the space of $\mathcal{F}_0$-measurable random variables with $\E|X_0|^{\tilde{p}}<\8$ for any $\tilde{p}>0$.\\
 \indent In this work, we study  McKean-Vlasov SDEs with common noise $W^0$ of the form
 \begin{equation}\label{eq1.2}
 	dX_t =b(X_t,{\cal{L}}_{X_t}^1)dt +\sigma(X_t,{\cal{L}}_{X_t}^1) dW_t + \sigma^0(X_t,{\cal{L}}_{X_t}^1) dW^0_t, \qquad X_0\in {L}_0^{\tilde{p}}(\mathbb{R}^d),
 \end{equation}
where $W^0$ is an $m^0$-dimensional Brownian motion and $({\cal L}_{X_t}^1)_{t\in[0,T]}$ denotes the flow of marginal conditional distributions of $X$ given the common noise.  In contrast to the McKean-Vlasov SDEs without common noise (\ref{eq1.1}), the marginal conditional distributions are no longer deterministic. The classical notion of propagation of chaos can be intuitively viewed as the idea that in a large network of $N$ interacting particles, these particles will gradually become independent as $N\rightarrow\8$. When all particles have a common random source, (i.e., all particles are influenced by the track of the common noise $W^0$,) it is reasonable that they are asymptotically independent. However, for the information generated by the common noise, it seems unlikely to become asymptotically independent as $N\rightarrow\8$. To put it another way, as $N\rightarrow\8$, it is expected that the empirical distribution of the particles will converge to the common conditional distribution of each particle given a common source of $W^0$. For more details on this topic, please refer to \cite{C18,CR18}.\\
\indent The simulation of McKean-Vlasov SDEs typically involves two steps. The first step is to use the empirical measure $$\mu_{t}^{{X},N}(dx):=\frac{1}{N}\sum_{j=1}^{N}\delta_{{X}_{t}^{j,N}}(dx),$$
to approximate the measure ${\cal L}_{X_t}^1$, where $\delta_{{X}_{t}^{j,N}}$ denotes the Dirac measure at point ${{X}_{t}^{j,N}}$ and $X^{i,N}$  is the solution to interacting particle system
\begin{equation}\label{eq1.3} 
	dX_t^{i,N} =b(X_t^{i,N},\mu_t^{X,N})dt + \sigma(X_t^{i,N},\mu_t^{X,N}) dW_t^i + \sigma^0(X_t^{i,N},\mu_t^{X,N}) dW^0_t,\quad X_0^{i,N}=X^i_0\in {L}_0^{\tilde{p}}(\mathbb{R}^d),
\end{equation}
where $W^i$ and $W^0$ are independent Brownian motions for any $i\in\{1,\dots, N\}.$ The second step is to construct a reasonable time-step scheme to discretise the particle system $(X^{i,N})_{\{i=1,\dots,N\}}$.\\
\indent The existence and uniqueness theory of strong solutions for McKean-Vlasov SDEs with coefficients exhibiting linear growth and satisfying Lipschitz-type conditions (both in the state and the measure components) is well-established (see\cite{S91,L18,Hu21}). In the case of the drift exhibiting super-linear growth while the diffusion maintains linear growth, it is known that McKean-Vlasov SDEs admit a unique strong solution\cite{R17}. The theories about well-posedness for McKean-Vlasov SDEs whose coefficients exhibit super-linear growth in the state variable have been extensively studied in various literature\cite{B20,Me20,K22}. For results of further advancements on the existence and uniqueness of weak and strong solutions, please refer to \cite{B18,M20,H21,Ka22,N24} and the references therein.\\
\indent In recent years, McKean-Vlasov SDEs have garnered substantial attention largely due to their wide-ranging applicability in diverse technical domains and mathematical models. The study of the McKean-Vlasov equations with common noise has become a focal point, as demonstrated in \cite{P16}, where a limiting equation for individual particles interacting with each other is presented. The significance of common noise in McKean-Vlasov SDE is further underscored by their applications in modeling complex systems. In the realm of financial systems, these models can effectively capture the dynamics of contagion and common exposure risks within large-scale networks, as demonstrated in \cite{L21}. Further literature, in which the motivation to consider the class McKean-Vlasov SDEs of a common noise source, plays a crucial role in McKean-Vlasov SDEs and associated interacting equations extends to the field of neuroscience. In this context, McKean-Vlasov equations are used to characterize the voltage fluctuations of a representative neuron in a complex network (see \cite{E23,U18}). The well-posedness theories of McKean-Vlasov SDEs with common noise have been established in \cite{CR18,K22,HW21}.\\
\indent Concerning classical SDEs, there is a large body of literature investigating numerical approximation\cite{B22,K92}. The scenario where coefficients exhibit super-linear growth has attracted significant attention and is now well-explored in the literature. Building on this, prior studies have introduced tamed EM schemes tailored for local Lipschitz drift coefficients that may potentially exhibit super-linear growth in \cite{H12,S13}. This approach was subsequently expanded to accommodate diffusion coefficients with super-linear growth, as detailed in references \cite{K19,H15}.\\
\indent The pioneering works on the numerical approximation of McKean-Vlasov SDEs within a continuous framework were presented in \cite{B97}. Conventional numerical approximation methods for McKean-Vlasov SDEs with Lipschitz continuous coefficients are well-documented, and strong convergence results with a convergence order of 1/2 have been established using various numerical schemes, e.g., explicit EM scheme \cite{L22}, tamed scheme \cite{K22}, projected scheme \cite{B18}, Milstein scheme \cite{BJ20,B21}. The investigation into stable time-stepping schemes for the tamed EM scheme and an implicit scheme for interacting particle systems (in the absence of common noise) was initiated in \cite{D22}. This study was conducted under the premise that the drift term is permitted to exhibit super-linear growth in the state component, while the diffusion term adheres to global Lipschitz continuity in both the state variable and the measure component.\\
\indent This paper is centered on a time-stepping scheme to approximate the particle system in equation (\ref{eq1.3}) by using an adaptive timestep. Fang utilized this approach to investigate various properties of SDEs in both finite and infinite time horizons \cite{W20}. After that, the methodology was introduced by Reisinger and Stockinger in their article, where they presented an adaptive time-stepping scheme for McKean-Vlasov SDEs with super-linear growth in the drift and diffusion, assuming only a monotonicity condition, (see \cite{C22}). Building upon this foundation, this paper extends the study to McKean-Vlasov SDEs with a common noise in both finite and infinite time. A key requirement of the paper is that both the drift term and the diffusion term grow superlinearly in state variables.

 In summary, under fairly general assumptions, this paper investigates the adaptive EM schemes for the McKean-Vlasov SDEs with common noise and shows the strong convergence rates in finite and infinite horizons.

The paper is structured as follows. In Section \ref{section2}, we introduce some basic notations and probabilistic framework which gives background results needed throughout this paper. In Section \ref{section3}, we describe the particle method and the adaptive scheme and give the main results both in  finite and infinite horizons. In Section \ref{section4}, we present several numerical examples to support our theoretical results.

\section{Mathematical preliminaries}\label{section2}\indent

$Notations.\quad$ Let $(\mathbb{R}_d,\langle \cdot,\cdot\rangle,|\cdot|)$ represent the $d$-dimensional Euclidean space. For all linear operators (e.g. matrices $A\in\RR^{d\times m}$) appearing in this article, we will use the standard Hilbert–Schmidt norm denoted by $\|\cdot\|$. The transpose of a matrix $A$ will be denoted by $A^\ast.$ In addition, we use ${\cal P}(\mathbb R^d)$ to denote the family of all probability measures on $(\mathbb R^d,{\cal B} (\mathbb R^d))$, where ${\cal B}(\mathbb R^d)$ denotes the Borel $\sigma$-field over $\mathbb R^d$, and define the subset of probability measures with finite $p$th moment by
\begin{equation*}
{\cal P}_p(\mathbb R^d):=\left\{\mu\in{\cal P}(\mathbb R^d):\left(\int_{\mathbb R^d}|x|^p\mu (dx)\right)^{1/(1\vee p)}<\8\right\}.
\end{equation*}
For any $\mu,\nu\in {\cal P}_p(\RR^d)$, the $L^p$-Wasserstein distance between $\mu$ and $\nu$ is defined as 
\begin{equation*}
	{\cal W}_p(\mu,\nu):=\inf_{\pi\in {\cal C}(\mu,\nu)}\left(\int_{\RR^d\times\RR^d}|x-y|^p\pi\left(dx,dy\right)\right)^{1/(1\vee p)},
\end{equation*}
where ${\cal C}(\mu,\nu)$ is the set of all the couplings of $\mu$ and $\nu$, i.e., $\pi\in C(\mu,\nu)$ if and only if $\pi(\cdot,\RR^d)=\mu(\cdot)$ and $\pi(\RR^d,\cdot)=\nu(\cdot)$.  ${\cal P}_p(\RR^d)$ is a Polish space under the $L^p$-Wasserstein metric. For a given $p\geq1,L_p^0(\RR^d)$ will denote the space of $\RR^d$-valued, ${\cal F}_0$-measurable random variables $X$ satisfying $\E|X|^p<\8.$ Furthermore, ${\cal S}^p([0, T])$ refers to the space of $\RR^d$-valued continuous, ${\cal F}$-adapted processes, defined on the interval $[0, T],$ for $T>0,$ with finite $p$-th moments.\\\\
\indent$Probabilistic\ framework.\quad$ We give the complete probability spaces $(\Omega^0,{\cal F}^0,P^0)$ and $(\Omega^1,{\cal F}^1,P^1),$ endowed with two right-continuous and complete filtrations $\mathbb{F}^0=({\cal F}_t^0)_{t\geq0}$ and $\mathbb{F}^1=({\cal F}_t^1)_{t\geq0}.$ Here, we assume the Brownian motion $W^0$ will be constructed on the $(\Omega^0,{\cal F}^0,P^0)$ and $W$ and $W^i$ will be constructed on the $(\Omega^1,{\cal F}^1,P^1).$ We then define the product space $(\Omega,{\cal F},P),$ where  $\Omega=\Omega^0\times\Omega^1,$ $({\cal F},P)$ is the completion of $({\cal F}^0\otimes{\cal F}^1,P^0\otimes P^1)$ and $\mathbb{F}=({\cal F}_t)_{t\geq0}$ is the complete and
right continuous augmentation of $({\cal F}_t^0\otimes{\cal F}_t^1)_{t\geq 0}.$ It is known from Lemma 2.4 in \cite{CR18} that for a given random variable $\Omega\ni X\rightarrow \RR^d$, equipped with the completion ${\cal F},$ the mapping ${\cal L}^1_X:\Omega^0\ni\omega^0\mapsto{\cal L}_{X(\omega^0,\cdot)},$ is $P^0$-almost surely well defined and forms a random variable from $(\Omega^0,{\cal F}^0,P^0)$ into ${\cal P}_2(\RR^d )$, ${\cal L}_X^1$ can also be seen as a conditional law of $X$ given $\mathcal{F}^0$.\\

\section{Adaptive time-stepping strategies}\label{section3}

\subsection{Definition of adaptive EM scheme}\indent

\indent In this subsection, we introduce the adaptive EM scheme for the interacting particle system (\ref{eq1.3}) that is associated with (\ref{eq1.2}), and we use a scheme that allows the computation of an empirical measure $\mu_t^{\hat{X},N}$ for all $t$. We start at $t=0$ with $\hat{X}_0^{i,N}=\xi$ for all $i\in \{1,\dots,N\}$. At step $n\geq0,$ now we perform the EM scheme given by
\begin{equation}
\hat{X}_{t_{n+1}}^{i,N}=\hat{X}_{t_n}^{i,N}+b\big(\hat{X}_{t_n}^{i,N},\mu_{t_n}^{\hat{X},N}\big)h_n^{\min}+\sigma\big(\hat{X}_{t_n}^{i,N},\mu_{t_n}^{\hat{X},N}\big)\Delta W_{t_n}^i+\sigma^0\big(\hat{X}_{t_n}^{i,N},\mu_{t_n}^{\hat{X},N}\big)\Delta W_{t_n}^0,
\end{equation}
where $\Delta W_{t_n}^i=W_{t_{n+1}}^i-W_{t_n}^i$ and $\Delta W_{t_n}^0=W_{t_{n+1}}^0-W_{t_n}^0,$ and the step $h_n^{\min}:=\min\{h_n^1,\dots,h_n^N\},$ for each particle the size of the adaptive time-step $h_n^i:=h(\hat{X}_{t_n}^{i,N},\mu_{t_n}^{\hat{X},N}),$ $t_{n+1}=t_n+h_n^{\min}$ and $t_n$ increases as $n$ increases until $n=M,$ such that $t_M\geq T,$ for each $i\in \{1,\dots,N\}$, where
\begin{equation}\label{3.2}
\mu_{t_n}^{\hat{X},N}(dx):=\frac{1}{N}\sum_{j=1}^{N}\delta_{\hat{X}_{t_n}^{j,N}}(dx).
\end{equation} 
We use the notation $\underline{t}:=\max\{t_n:t_n\leq t\}$ to represent the nearest time point before time $t$ and $n_t:=\max\{n:t_n\leq t\}$ to denote the number of step approximations up to time $t$. Then we introduce the piecewise constant interpolant process $\bar{X}_t^{i,N}:=\hat {X}_{\underline{t}}^{i,N}$ and measure $\mu_{t}^{\bar{X},N}:=\mu_{\underline{t}}^{\hat{X},N}$, and also define the continuous interpolant process 
\begin{equation}\label{3.3}
\hat{X}_{t}^{i,N}=\bar{X}_{{t}}^{i,N}+b(\hat{X}_{\underline{t}}^{i,N},\mu_{\underline{t}}^{\hat{X},N})(t-\underline{t})+\sigma(\hat{X}_{\underline{t}}^{i,N},\mu_{\underline{t}}^{\hat{X},N})(W_{t}^i-W_{\underline{t}}^i)+\sigma^0(\hat{X}_{\underline{t}}^{i,N},\mu_{\underline{t}}^{\hat{X},N})(W_t^0-W_{\underline{t}}^0),
\end{equation}
so that $\big(\hat{X}_{t}^{i,N}\big)_{t\in [0,T]}$ is the solution to the SDE
\begin{align}
d\hat{X}_{t}^{i,N}&=b(\hat{X}_{\underline{t}}^{i,N},\mu_{\underline{t}}^{\hat{X},N})dt+\sigma(\hat{X}_{\underline{t}}^{i,N},\mu_{\underline{t}}^{\hat{X},N})dW_{t}^i+\sigma^0(\hat{X}_{\underline{t}}^{i,N},\mu_{\underline{t}}^{\hat{X},N})dW_t^0\nonumber\\
&=b(\bar{X}_t^{i,N},\mu_{t}^{\bar{X},N})dt+\sigma(\bar{X}_t^{i,N},\mu_{t}^{\bar{X},N})dW_{t}^i+\sigma^0(\bar{X}_t^{i,N},\mu_{t}^{\bar{X},N})dW_t^0.
\end{align}

\subsection{Boundedness and convergence in finite horizon}\indent
\indent This subsection investigates the boundedness and convergence of adaptive scheme in finite horizon. 
\begin{assp}\label{assp1}The mappings $ b:\mathbb R^d \times{\cal P}_2 (\mathbb R^d)\rightarrow\mathbb R^d$, $ \sigma:\mathbb R^d \times{\cal P}_2 (\mathbb R^d)\rightarrow\mathbb R^{d\times m}$ and $ \sigma^0:\mathbb R^d \times{\cal P}_2 (\mathbb R^d)\rightarrow\mathbb R^{d\times m_0}$ are measurable and satisfy:\\\\
(1) For some $p\geq2$, there exists a constant $L>0$ such that
\begin{align*}
&\langle x-x{'},b(x,\mu)-b(x{'},\mu{'})\rangle+(p-1)\|\sigma(x,\mu)-\sigma(x{'},\mu{'}) \|^2\\
&+(p-1) \|\sigma^0(x,\mu)-\sigma^0(x{'},\mu{'}) \|^2\le L\big (| x - x{'} |^2+\mathcal{W}_2^2(\mu,\mu{'})\big),
\end{align*} 
for any $x,x{'}\in \mathbb R^d and\ \mu,\mu{'}\in {\cal P}_2(\mathbb R^d)$.\\
(2) There exist constants $L>0,q>0$ such that
\begin{equation*}
|b(x,\mu)-b(x{'},\mu{'})|\le L\big[(1+|x|^q+|x{'}|^q)|x-x{'}|+\mathcal{W}_2(\mu,\mu{'})\big],
\end{equation*}for any $x,x{'}\in \mathbb R^d and\ \mu,\mu{'}\in {\cal P}_2(\mathbb R^d).$
\end{assp}

\begin{rmk}\label{rmk1}
	Due to Assumption \ref{assp1}, there exists a constant $K:=K(L)>0$ such that
	\begin{align*}
	\|\sigma(x,\mu)-\sigma(x{'},\mu{'}) \|^2+ \|\sigma^0(x,\mu)-\sigma^0(x{'},\mu{'}) \|^2\le K\big[(1+|x|^{q}+|x{'}|^{q})| x - x{'} |^2+\mathcal{W}_2^2(\mu,\mu{'})\big],
	\end{align*}
for any  $t\in[0,T],x,x{'}\in \mathbb R^d$ and $\mu,\mu{'}\in {\cal P}_2(\mathbb R^d).$ And
\begin{align*}
		\big\langle& x,b(x,\mu)\rangle+\frac{(p-1)}{2}\|\sigma(x,\mu)\|^2+\frac{(p-1)}{2}\|\sigma^0(x,\mu)\|^2\\
		\leq&\big\langle x,b(x,\mu)-b(0,\delta_0)\big\rangle+\big\langle x,b(0,\delta_0)\big\rangle+(p-1)\|\sigma(x,\mu)-\sigma(0,\delta_0)\|^2+(p-1)\|\sigma(0,\delta_0)\|^2\nonumber\\
		&+(p-1)\|\sigma^0(x,\mu)-\sigma^0(0,\delta_0)\|^2+(p-1)\|\sigma^0(0,\delta_0)\|^2\nonumber\\
		\leq&L|x|^2+L\mathcal{W}_2^2(\mu,\delta_0)+|x||b(0,\delta_0)|+(p-1)\|\sigma(0,\delta_0)\|^2+(p-1)\|\sigma^0(0,\delta_0)\|^2\nonumber\\
		\leq&L|x|^2+L\mathcal{W}_2^2(\mu,\delta_0)+\varepsilon|x|^2+\frac{1}{\varepsilon}|b(0,\delta_0)|^2+(p-1)\|\sigma(0,\delta_0)\|^2+(p-1)\|\sigma^0(0,\delta_0)\|^2\nonumber\\
		=&(\varepsilon+L)|x|^2+L\mathcal{W}_2^2(\mu,\delta_0)+\frac{1}{\varepsilon}|b(0,\delta_0)|^2+(p-1)\|\sigma(0,\delta_0)\|^2+(p-1)\|\sigma^0(0,\delta_0)\|^2\\
		\leq& K_1\big(1+|x|^2+{\cal W}_2^2(\mu,\delta_0)\big),
\end{align*}
	for any $x\in \mathbb R^d $ and $ \mu\in {\cal P}_2(\mathbb R^d),$ where constant $K_1:=(\varepsilon+L)\vee (\frac{1}{\varepsilon}|b(0,\delta_0)|^2+(p-1)\|\sigma(0,\delta_0)\|^2+(p-1)\|\sigma^0(0,\delta_0)\|^2).$ 
\end{rmk}\indent
\indent In the preceding discussion, a set of assumptions has been posited to ensure that the well-posedness of the equation under consideration remains intact. Furthermore, the establishment of moment bounds and the demonstration of strong convergence for the adaptive EM scheme are pivotal objectives that will be addressed in the forthcoming sections.\\
\indent Compared to the classic EM scheme, which employs a constant time step size $h$, the present analysis incorporates an adaptive scheme wherein the timestep $h(x,\mu)$ is contingent upon the state variable $x$. This necessitates an adjustment to the timestep management strategy. We will consider a timestep function $h^{\delta}(x,\mu)$ that is governed by a tunable parameter $\delta$ with $0<\delta\leq1$. Subsequently, the convergence analysis hinges solely on the examination of the limit as $\delta\rightarrow0$.
\begin{assp}\label{assp2}(Adaptive timestep)\\
(1) The adaptive timestep function $h:\RR^d \rightarrow \RR^+$ is continuous and positive, and there exists a constant $C_h>0$ such that
\begin{equation*}
	h(x,\mu)= \frac{C_h}{\big(1+|b(x,\mu)|\|\sigma(x,\mu)\|+|b(x,\mu)|\|\sigma^0(x,\mu)\|+|x|^q\big)^{2}},
\end{equation*}
where $q$ is defined in Assumption \ref{assp1}.\\
(2) The timestep function $h^\delta:\RR^d \rightarrow \RR^+,0<\delta \leq1$ satisfies
\begin{equation*}
	\delta \min(T,h(x,\mu))\leq h^\delta(x,\mu)\leq \min (\delta T,h(x,\mu)).
\end{equation*}
(3)	There exist constants $\alpha_1,\alpha_2,\beta,\varpi>0$ such that 
\begin{equation}\label{C3.4}
	h(x,\mu)\geq\big(\alpha_1|x|^{\varpi}+\alpha_2\mathcal{W}_2^{\varpi}(\mu,\delta_0)+\beta\big)^{-1},\quad\forall x\in\RR^d.
\end{equation}
\end{assp}
\begin{rmk}\label{rmk3.2}
	According to Assumption \ref{assp1}(2) and Remark \ref{rmk1}, there exist constants $C_1, C_2,C_3,C_4>0$ such that
	\begin{equation*}
		|b(x,\mu)|\le C_1\big[(1+|x|^{q+1})+\mathcal{W}_2(\mu,\delta_0)\big],
	\end{equation*}
	\begin{equation*}
		\|\sigma(x,\mu)\|^2+ \|\sigma^0(x,\mu)\|^2\le C_2\big[(1+|x|^{q+2})+\mathcal{W}_2^2(\mu,\delta_0)\big],
	\end{equation*} for any $ x\in \mathbb R^d and\ \mu\in {\cal P}_2(\mathbb R^d).$
The constraints on the timestep function $h$ in Assumption \ref{assp2} could ensure that
	\begin{equation*}
		|b(x,\mu)|^2h(x,\mu)\leq C_3(1+|x|^2+\mathcal{W}_2^2(\mu,\delta_0)),
	\end{equation*}
	\begin{equation}\label{C3.5}
		|b(x,\mu)|\|\sigma(x,\mu)\|h(x,\mu)^{\frac{1}{2}}+|b(x,\mu)|\|\sigma^0(x,\mu)\|h(x,\mu)^{\frac{1}{2}}\leq C_4,
	\end{equation}
which will play an important role in proving the  moment boundedness of numerical solution in Theorem \ref{theorem1}.
\end{rmk}
\begin{rmk}
	We emphasize that the corresponding time-step function $h(x,\mu)$ can be chosen under the  Assumption \ref{assp2}. For example, let $b(x,\mu)=-2x^3+\int_{\mathbb{R}}^{}x\mu(dx),$ $\sigma(x,\mu)=0.5x^2+\int_{\mathbb{R}}^{}x\mu(dx),$ $\sigma^0(x,\mu)=0.5x^2+\int_{\mathbb{R}}^{}x\mu(dx)$. Due to (\ref{C3.5}), then $h(x,\mu)=(1+|x|^5+\int_{\mathbb{R}}^{}|x|^2\mu(dx))^{-2}$ and (\ref{C3.4}) is satisfied by taking $C_h=1,\alpha_1=\alpha_2=1,\beta=2,\varpi=10$. Then more details can be seen in Section 4.
\end{rmk}
The proof of the following lemma can be found  in Theorem 2.1 of \cite{K22}.
\begin{lemma}\label{lemma1}
	Let Assumption \ref{assp1} hold. Then for all $0<r\leq{\tilde{p}}$, (\ref{eq1.2}) admits a unique strong solution ${X}_t$ such that $$\sup_{t\in[0,T]}E\left[\big|{X}_t\big|^r\right]\leq C_{X_0}.$$
\end{lemma}
\begin{theorem}\label{theorem1}
Let Assumption \ref{assp1} and Assumption \ref{assp2} hold. Then $T$ is almost surely attainable, and for all $0<p<\tilde{p}$, there exists a constant $C>0$, which depends on $T$ and $p$, such that 
$$\max_{i\in \{1,\dots,N\}}\sup_{t\in[0,T]}E\left[\big|\hat{X}_t^{i,N}\big|^p\right]\leq C.$$
\end{theorem}\indent
\indent\proof
Since the step size is variable, we need to show for any $T>0$ and almost all $\omega\in\Omega$, there exists finite $M(\omega)$ such that $t_{M(\omega)}\geq T$, i.e. $\mathbb{P}\{\exists M(\omega)<\8,\ s.t.\ t_{M(\omega)\geq T}\}=1$.
Then, define a stopping time $\tau_R:=\inf\{t\geq0:\max_{i\in\{1,\dots,N\}}\big|\hat{X}_t^{i,N}\big|\geq R\}$ for each $R>0.$ By It\^{o}'s formula, we have
\begin{align*}
	\big(1+\big|\hat{X}_{t\wedge\tau_R}^{i,N}\big|^2\big)^{p/2}=&\big(1+\big|\hat{X}_{0}^{i,N}\big|^2\big)^{p/2}+p\int_{0}^{t\wedge\tau_R}\big(1+\big|\hat{X}_{s}^{i,N}\big|^2\big)^{p/2-1}\big\langle\hat{X}_{s}^{i,N},b(\bar{X}_{s}^{i,N},\mu_{s}^{\bar{X},N})\big\rangle ds\\
	&+p\int_{0}^{t\wedge\tau_R}\big(1+\big|\hat{X}_{s}^{i,N}\big|^2\big)^{p/2-1}\big\langle\hat{X}_{s}^{i,N},\sigma(\bar{X}_{s}^{i,N},\mu_{s}^{\bar{X},N})dW_s^i\big\rangle\\
	&+p\int_{0}^{t\wedge\tau_R}\big(1+\big|\hat{X}_{s}^{i,N}\big|^2\big)^{p/2-1}\big\langle\hat{X}_{s}^{i,N},\sigma^0(\bar{X}_{s}^{i,N},\mu_{s}^{\bar{X},N})dW_s^0\big\rangle\\
	&+\frac{p(p-2)}{2}\int_{0}^{t\wedge\tau_R}\big(1+\big|\hat{X}_{s}^{i,N}\big|^2\big)^{p/2-2}\big|\sigma^\ast(\bar{X}_{s}^{i,N},\mu_{s}^{\bar{X},N})\hat{X}_{s}^{i,N}\big|^2ds\\
	&+\frac{p(p-2)}{2}\int_{0}^{t\wedge\tau_R}\big(1+\big|\hat{X}_{s}^{i,N}\big|^2\big)^{p/2-2}\big|\sigma^{0,\ast}(\bar{X}_{s}^{i,N},\mu_{s}^{\bar{X},N})\hat{X}_{s}^{i,N}\big|^2ds\\
	&+\frac{p}{2}\int_{0}^{t\wedge\tau_R}\big(1+\big|\hat{X}_{s}^{i,N}\big|^2\big)^{p/2-1}\big\|\sigma(\bar{X}_{s}^{i,N},\mu_{s}^{\bar{X},N})\big\|^2ds\\
	&+\frac{p}{2}\int_{0}^{t\wedge\tau_R}\big(1+\big|\hat{X}_{s}^{i,N}\big|^2\big)^{p/2-1}\big\|\sigma^{0}(\bar{X}_{s}^{i,N},\mu_{s}^{\bar{X},N})\big\|^2ds,
\end{align*}
almost surely for any $t\in[0,T],i\in\{1,\dots,N\}$. Consequently, by taking the expectation and using Remark \ref{rmk1},  Young's inequality, we deduce
\begin{align}\label{3.5}
	\E\left[\big(1+\big|\hat{X}_{t\wedge\tau_R}^{i,N}\big|^2\big)^{p/2}\right]\leq&\E\left[\big(1+\big|\hat{X}_{0}^{i,N}\big|^2\big)^{p/2}\right]\nonumber\\
	&+p\E\Big[\int_{0}^{t\wedge\tau_R}\big(1+\big|\hat{X}_{s}^{i,N}\big|^2\big)^{p/2-1}\Big(\big\langle\bar{X}_{s}^{i,N},b(\bar{X}_{s}^{i,N},\mu_{s}^{\bar{X},N})\big\rangle\nonumber\\
	&+\frac{(p-1)}{2}\|\sigma(\bar{X}_{s}^{i,N},\mu_{s}^{\bar{X},N})\|^2+\frac{(p-1)}{2}\|\sigma^0(\bar{X}_{s}^{i,N},\mu_{s}^{\bar{X},N})\|^2\Big)ds\Big]\nonumber\\
	&+p\E\Big[\int_{0}^{t\wedge\tau_R}\big(1+\big|\hat{X}_{s}^{i,N}\big|^2\big)^{p/2-1}\big\langle\hat{X}_s^{i,N}-\bar{X}_{s}^{i,N},b(\bar{X}_{s}^{i,N},\mu_{s}^{\bar{X},N})\big\rangle ds\Big]\nonumber\\
	\leq&\E\left[\big(1+\big|\hat{X}_{0}^{i,N}\big|^2\big)^{p/2}\right]+C_p\E\Big[\int_{0}^{t\wedge\tau_R}\big(1+\big|\hat{X}_{s}^{i,N}\big|^2\big)^{p/2}ds\Big]\nonumber\\
	&+C_p\E\Big[\int_{0}^{t\wedge\tau_R}\big(1+\big|\hat{X}_{s}^{i,N}\big|^2\big)^{p/2-1}\big(1+\big|\bar{X}_{s}^{i,N}\big|^2\big)ds\Big]\nonumber\\
	&+C_p\E\Big[\int_{0}^{t\wedge\tau_R}\big(1+\big|\hat{X}_{s}^{i,N}\big|^2\big)^{p/2-1}{\cal W}_2^2(\mu_{s}^{\bar{X},N},\delta_0)ds\Big]\nonumber\\
	&+C_p\E\left[\int_{0}^{t\wedge\tau_R}\big|\hat{X}_s^{i,N}-\bar{X}_{s}^{i,N}\big|^{p/2}\big|b(\bar{X}_{s}^{i,N},\mu_{s}^{\bar{X},N})\big|^{p/2}ds\right].
\end{align}
Hence, by  (\ref{3.3}) and Assumption \ref{assp2}, taking $\pi={s\wedge\tau_R}$ gives 
\begin{align*}
	\E\Big[&\big|\hat{X}_\pi^{i,N}-\bar{X}_\pi^{i,N}\big|^{p/2}\big|b(\bar{X}_\pi^{i,N},\mu_\pi^{\bar{X},N})\big|^{p/2}\Big]\\
	\leq&C_p\E\left[\big|b(\bar{X}_\pi^{i,N},\mu_\pi^{\bar{X},N})\big|^{p}(\pi-\underline{\pi})^{p/2}\right]\\
	&+C_p\E\left[\big|\sigma(\bar{X}_{\pi}^{i,N},\mu_{\pi}^{\bar{X},N})(W_\pi^i-W_{\underline{\pi}}^i)\big|^{p/2}\big|b(\bar{X}_{\pi}^{i,N},\mu_{\pi}^{\bar{X},N})\big|^{p/2}\right]\\
	&+C_p\E\left[\big|\sigma^0(\bar{X}_{\pi}^{i,N},\mu_{\pi}^{\bar{X},N})(W_\pi^0-W_{\underline{\pi}}^0)\big|^{p/2}\big|b(\bar{X}_{\pi}^{i,N},\mu_{\pi}^{\bar{X},N})\big|^{p/2}\right]\\
	\leq&C_p\E\left[\big(1+\big|\hat{X}_{\underline{\pi}}^{i,N}\big|^2+{\cal W}_2^2(\mu_{\pi}^{\bar{X},N},\delta_0)\big)^{p/2}\right]\\
	&+C_p\E\left[\E\left[\big|\sigma(\bar{X}_{\pi}^{i,N},\mu_{\pi}^{\bar{X},N})(W_\pi^i-W_{\underline{\pi}}^i)\big|^{p/2}\big|b(\bar{X}_{\pi}^{i,N},\mu_{\pi}^{\bar{X},N})\big|^{p/2}\Big|{\cal F}_{\underline{\pi}}\right]\right]\\
	&+C_p\E\left[\E\left[\big|\sigma^0(\bar{X}_{\pi}^{i,N},\mu_{\pi}^{\bar{X},N})(W_\pi^0-W_{\underline{\pi}}^0)\big|^{p/2}\big|b(\bar{X}_{\pi}^{i,N},\mu_{\pi}^{\bar{X},N})\big|^{p/2}\Big|{\cal F}_{\underline{\pi}}\right]\right]\\
	\leq&C_p\E\left[\big(1+\big|\bar{X}_{\pi}^{i,N}\big|^2+{\cal W}_2^2(\mu_{\pi}^{\bar{X},N},\delta_0)\big)^{p/2}\right],
\end{align*}
where $\E\left[\big|W_\pi^i-W_{\underline{\pi}}^i\big|^p\Big|{\cal F}_{\underline{\pi}}\right]\leq C_p(\pi-\underline{\pi})^{p/2},$ and $\E\Big[{\cal W}_2^2(\mu_{\pi}^{\bar{X},N},\delta_0)\Big]=\E\Big[\big|\bar{X}_\pi^{i,N}\big|^2\Big]$.
By substituting the aforementioned equations into equation (\ref{3.5}) and applying Young's inequality, we arrive at
\begin{align*}
		\sup_{s\in[0,t]}\E\left[\big(1+\big|\hat{X}_{s\wedge\tau_R}^{i,N}\big|^2\big)^{p/2}\right]\leq\E\left[\big(1+\big|\hat{X}_{0}^{i,N}\big|^2\big)^{p/2}\right]+C_p\int_{0}^{t}\sup_{r\in[0,s]}\E\Big[\big(1+\big|\hat{X}_{r\wedge\tau_R}^{i,N}\big|^2\big)^{p/2}\Big]ds.
\end{align*}
Then using Gronwall's inequality, we get
\begin{align*}
	\sup_{t\in[0,T]}\E\Big[|\hat{X}_{t\wedge\tau_R}^{i,N}|^p\Big]\leq C_{p,t}.
\end{align*}
Note that
\begin{align*}
	\PP(\tau_R\leq T)&\leq\sum_{i=1}^{N}\PP\left(\big|\hat{X}_{T\wedge\tau_R}^{i,N}\big|>R\right)\\
	&\leq N\PP\Big(\max_{i\in\{1,\dots,N\}}|\hat{X}_{T\wedge\tau_R}^{i,N}|>R\Big)\\
	&\leq \frac{N}{R^2}\E\left[\max_{i\in\{1,\dots,N\}}|\hat{X}_{T\wedge\tau_R}^{i,N}|^2\right]\leq\frac{CN}{R^2}.
\end{align*}
Hence
\begin{align*}
	\PP\left(\max_{i\in\{1,\dots,N\}}\sup_{t\in[0,T]}|\hat{X}_{t}^{i,N}|<R\right)=1-\PP(\tau_R\leq T)\rightarrow1,\:as\ R\rightarrow\8.
\end{align*}
Therefore, $\tau_R\uparrow\8$ as $R\uparrow\8$, $\max_{i\in\{1,\dots,N\}}\sup_{t\in[0,T]}|\hat{X}_{t}^{i,N}|<\8$ and $T$ is attainable. Then using the stability up to time $T\wedge\tau_R,$ the claim follows from Fatou's lemma.
\eproof

Given time $t$, the step-size depends on $\hat{X}_t^{i,N}$ and $\mu_{t}^{\hat{X},N}$, thus for each particle step size number $M_T^i$ until time $T$ is a random variable. Consequently the total number of steps $\E\left[\frac{1}{N}\sum_{i=1}^{N}M_T^i\right]$ makes the following estimate.

 \begin{prop}
 	Let Assumption \ref{assp1} and Assumption \ref{assp2} hold. Then  there exists $C>0$ independent of $\delta$ and $N$ such that $$\E\left[\frac{1}{N}\sum_{i=1}^{N}M_T^i\right]\leq C\delta^{-1}.$$
\end{prop}\indent
	\indent\proof
	Obviously, we see
 	\begin{align*}
 		M_T^i&\leq 1+T\sup_{0\leq t\leq T}\frac{1}{h^\delta(\hat{X}_t^{i,N},\mu_{t}^{\hat{X},N})}\\
 		&\leq 1+T\delta^{-1}\left(\sup_{0\leq t\leq T}\max\big(h^{-1}(\hat{X}_t^{i,N},\mu_{t}^{\hat{X},N}),T^{-1}\big)\right)\\
 		&\leq T\delta^{-1}\left(\alpha_1\sup_{0\leq t\leq T}\big|\hat{X}_t^{i,N}\big|^\varpi+\alpha_2{\cal W}_2^\varpi(\mu_{t}^{\hat{X},N},\delta_0)+\beta+(1+\delta)T^{-1}\right).
 	\end{align*}
 We can get the result immediately by Theorem \ref{theorem1}.
 \eproof

\indent To demonstrate the approximation of the particle system, we now introduce the propagation of chaos result. Initially, since we have already defined the interacting particle system in (\ref{eq1.3}), we shall focus exclusively on the system of noninteracting particles:
\begin{equation}\label{3.6}
	dX_t^i =b(X_t^i,{\cal{L}}_{X_t^i}^1)dt +\sigma(X_t^i,{\cal{L}}_{X_t^i}^1)dW_t^i +\sigma^0(X_t^i,{\cal{L}}_{X_t^i}^1)dW^0_t, 
\end{equation}
almost surely for any $t\in[0,T]$ and $i\in\{1,\dots ,N\}.$ Moreover, by Proposition 2.11 in \cite{CR18}, $$P^0\left[{\cal L}^1(X_t^i)={\cal L}^1(X_t^1) \ for \ all \ t\in[0,T]\right]=1.$$

The following result can be shown in the same way as shown in Theorem 3.2 in \cite{G24}, so we omit it here.
\begin{prop}\label{prop2}
	Let Assumption \ref{assp1} be satisfied with $p\in[2,r)$, then 
	\begin{equation*}
		\sup_{i\in\{1,\dots,N\}}\sup_{0\leq t\leq T}\E\Big[\big|X_t^i-X_t^{i,N}\big|^p\Big]\le C_{r,p,d}\left\{
		\begin{array}{l}
			N^{-1/2}+N^{-(r-p)/r},\quad\quad \quad\quad\quad if\  p>d/2\ and\ r\neq2p,\\
			N^{-1/2}\log(1+N)+N^{-(r-p)/r},\: if\  p=d/2\ and\ r\neq2p,\\
			N^{-p/d}+N^{-(r-p)/r},\quad   if\  p\in[2,d/2)\ and\ r\neq d/(d-p),
		\end{array}\right. 
	\end{equation*}
	where the constant $C_{r,p,d}>0$ does not depend on $N$.
\end{prop}

 \begin{theorem}(strong convergence)\label{theorem3}
 Let Assumption \ref{assp1} and Assumption \ref{assp2} hold. Then, for any $2\leq p<\tilde{p}$, there exists $C_{p,T}>0$ such that $$\max_{i\in\{1,\dots,N\}}\sup_{0\leq t\leq T}\E\Big[\big|\hat{X}_t^{i,N}-{X}_t^{i,N}\big|^p\Big]\leq C_{p,T}\delta^{p/2}.$$
 \end{theorem}\indent
\indent \proof
 First, define $e_t^i=\hat{X}_t^{i,N}-X_t^{i,N}$. We get 
 \begin{align*}
 	de_t^i=&\Big(b(\bar{X}_t^{i,N},\mu_t^{\bar{X},N})-b({X}_t^{i,N},\mu_t^{{X},N})\Big)dt+\Big(\sigma(\bar{X}_t^{i,N},\mu_t^{\bar{X},N})-\sigma({X}_t^{i,N},\mu_t^{{X},N})\Big)dW_t^i\\
 	&+\Big(\sigma^0(\bar{X}_t^{i,N},\mu_t^{\bar{X},N})-\sigma^0({X}_t^{i,N},\mu_t^{{X},N})\Big)dW_t^0.
 \end{align*}
 Subsequently, employing It\^{o}'s formula and the fact that $e_0^i=0,$ we derive 
 \begin{align*}
 	\big|e_t^i\big|^p\leq&p\int_{0}^{t}|e_s^i|^{p-2}\Big\langle e_s^i,b(\bar{X}_s^{i,N},\mu_s^{\bar{X},N})-b({X}_s^{i,N},\mu_s^{{X},N})\Big\rangle ds\\
 	&+\frac{p(p-1)}{2}\int_{0}^{t}|e_s^i|^{p-2}\big\|\sigma(\bar{X}_s^{i,N},\mu_s^{\bar{X},N})-\sigma({X}_s^{i,N},\mu_s^{{X},N})\big\|^2ds\\
 	&+\frac{p(p-1)}{2}\int_{0}^{t}|e_s^i|^{p-2}\big\|\sigma^0(\bar{X}_s^{i,N},\mu_s^{\bar{X},N})-\sigma^0({X}_s^{i,N},\mu_s^{{X},N})\big\|^2ds\\
 	&+p\int_{0}^{t}|e_s^i|^{p-2}\Big\langle e_s^i,\big(\sigma(\bar{X}_s^{i,N},\mu_s^{\bar{X},N})-\sigma({X}_s^{i,N},\mu_s^{{X},N})\big)dW_s^i\Big\rangle \\
 	&+p\int_{0}^{t}|e_s^i|^{p-2}\Big\langle e_s^i,\big(\sigma^0(\bar{X}_s^{i,N},\mu_s^{\bar{X},N})-\sigma^0({X}_s^{i,N},\mu_s^{{X},N})\big)dW_s^0\Big\rangle .
 \end{align*}
 Note that
 \begin{align*}
 	&\Big\langle e_s^i, b(\bar{X}_s^{i,N},\mu_s^{\bar{X},N})-b({X}_s^{i,N},\mu_s^{{X},N})\Big\rangle\\
 	&=\Big\langle e_s^i, b(\bar{X}_s^{i,N},\mu_s^{\bar{X},N})-b(\hat{X}_s^{i,N},\mu_s^{\hat{X},N})\Big\rangle+\Big\langle e_s^i, b(\hat{X}_s^{i,N},\mu_s^{\hat{X},N})-b({X}_s^{i,N},\mu_s^{{X},N})\Big\rangle.
 \end{align*}
Using the polynomial growth condition in Assumption \ref{assp1}, we know
 \begin{align}\label{3.7}
 	\Big\langle& e_s^i, b(\bar{X}_s^{i,N},\mu_s^{\bar{X},N})-b(\hat{X}_s^{i,N},\mu_s^{\hat{X},N})\Big\rangle\nonumber\\
 	&\leq \big|e_s^i\big| \left\{\tilde{Q}_L(\hat{X}_s^{i,N},\bar{X}_s^{i,N})\big|\hat{X}_s^{i,N}-\bar{X}_s^{i,N}\big|+{\cal W}_2(\mu_s^{\hat{X},N},\mu_s^{\bar{X},N})\right\}\nonumber\\
 	&\leq \frac{1}{2}\big|e_s^i\big|^2+\tilde{Q}_L(\hat{X}_s^{i,N},\bar{X}_s^{i,N})^2\big|\hat{X}_s^{i,N}-\bar{X}_s^{i,N}\big|^2+{\cal W}_2^2(\mu_s^{\hat{X},N},\mu_s^{\bar{X},N}),
 	\end{align}
 where $\tilde{Q}_L(\hat{X}_s^{i,N},\bar{X}_s^{i,N}):=L(1+\big|\hat{X}_s^{i,N}\big|^q+\big|\bar{X}_s^{i,N}\big|^q).$ 
 Similarly, 
 \begin{align*}
 	\big\|&\sigma(\bar{X}_s^{i,N},\mu_s^{\bar{X},N})-\sigma({X}_s^{i,N},\mu_s^{{X},N})\big\|^2+\big\|\sigma^0(\bar{X}_s^{i,N},\mu_s^{\bar{X},N})-\sigma^0({X}_s^{i,N},\mu_s^{{X},N})\big\|^2\\
 	&\leq2\big\|\sigma(\bar{X}_s^{i,N},\mu_s^{\bar{X},N})-\sigma(\hat{X}_s^{i,N},\mu_s^{\hat{X},N})\big\|^2+2\big\|\sigma(\hat{X}_s^{i,N},\mu_s^{\hat{X},N})-\sigma({X}_s^{i,N},\mu_s^{{X},N})\big\|^2\\
 	&\quad+2\big\|\sigma^0(\bar{X}_s^{i,N},\mu_s^{\bar{X},N})-\sigma^0(\hat{X}_s^{i,N},\mu_s^{\hat{X},N})\big\|^2+2\big\|\sigma^0(\hat{X}_s^{i,N},\mu_s^{\hat{X},N})-\sigma^0({X}_s^{i,N},\mu_s^{{X},N})\big\|^2.
 \end{align*}
 Combining Remark \ref{rmk1} on $\sigma$ and $\sigma^0$, we are able to deduce that
 \begin{align}\label{3.8}
 	\big\|&\sigma(\bar{X}_s^{i,N},\mu_s^{\bar{X},N})-\sigma(\hat{X}_s^{i,N},\mu_s^{\hat{X},N})\big\|^2+\big\|\sigma^0(\bar{X}_s^{i,N},\mu_s^{\bar{X},N})-\sigma^0(\hat{X}_s^{i,N},\mu_s^{\hat{X},N})\big\|^2\nonumber\\
 	&\leq2\tilde{Q}_K(\hat{X}_s^{i,N},\bar{X}_s^{i,N})^2\big|\hat{X}_s^{i,N}-\bar{X}_s^{i,N}\big|^2+2{\cal W}_2^2(\mu_s^{\hat{X},N},\mu_s^{\bar{X},N}),
 \end{align}
where $\tilde{Q}_K(\hat{X}_s^{i,N},\bar{X}_s^{i,N}):=K(1+\big|\hat{X}_s^{i,N}\big|^{q}+\big|\bar{X}_s^{i,N}\big|^{q}).$
According to Assumption \ref{assp1}, we know
\begin{align*}
	&\Big\langle e_s^i,b(\hat{X}_s^{i,N},\mu_s^{\hat{X},N})-b({X}_s^{i,N},\mu_s^{{X},N})\Big\rangle
	+(p-1)\big\|\sigma(\hat{X}_s^{i,N},\mu_s^{\hat{X},N})-\sigma({X}_s^{i,N},\mu_s^{{X},N})\big\|^2\\
	&+(p-1)\big\|\sigma^0(\hat{X}_s^{i,N},\mu_s^{\hat{X},N})-\sigma^0({X}_s^{i,N},\mu_s^{{X},N})\big\|^2\\
	&\leq L\big(|e_s^i|^{2}+{\cal W}_2^2(\mu_s^{\hat{X},N},\mu_s^{{X},N})\big).
\end{align*}
Using Young's inequality and all the estimates above, we  conclude that
\begin{align}\label{3.9}
	\big|e_t^i\big|^p\leq&Lp\int_{0}^{t}|e_s^i|^{p-2}\Big(|e_s^i|^{2}+{\cal W}_2^2(\mu_s^{\hat{X},N},\mu_s^{{X},N})\Big)ds+\frac{p}{2}\int_{0}^{t}|e_s^i|^{p}ds\nonumber\\
	&+\int_{0}^{t}|e_s^i|^{p-2}\big(\tilde{Q}_L(\hat{X}_s^{i,N},\bar{X}_s^{i,N})^2\big|\hat{X}_s^{i,N}-\bar{X}_s^{i,N}\big|^2+{\cal W}_2^2(\mu_s^{\hat{X},N},\mu_s^{\bar{X},N})\big)ds\nonumber\\
	&+2(p-1)\int_{0}^{t}|e_s^i|^{p-2}\big(\tilde{Q}_K(\hat{X}_s^{i,N},\bar{X}_s^{i,N})^2\big|\hat{X}_s^{i,N}-\bar{X}_s^{i,N}\big|^2+{\cal W}_2^2(\mu_s^{\hat{X},N},\mu_s^{\bar{X},N})\big)ds\nonumber\\
	&+p\int_{0}^{t}|e_s^i|^{p-2}\Big\langle e_s^i,\big(\sigma(\bar{X}_s^{i,N},\mu_s^{\bar{X},N})-\sigma({X}_s^{i,N},\mu_s^{{X},N})\big)dW_s^i\Big\rangle\nonumber \\
	&+p\int_{0}^{t}|e_s^i|^{p-2}\Big\langle e_s^i,\big(\sigma^0(\bar{X}_s^{i,N},\mu_s^{\bar{X},N})-\sigma^0({X}_s^{i,N},\mu_s^{{X},N})\big)dW_s^0\Big\rangle .
\end{align}
 Additionally, the definition of Wasserstein distance yields $$\E\left[{\cal W}_2^2(\mu_s^{\hat{X},N},\mu_s^{\bar{X},N})\right]\leq \E\left[\big|\hat{X}_s^{i,N}-\bar{X}_s^{i,N}\big|^2\right],\quad\E\left[{\cal W}_2^2(\mu_s^{\hat{X},N},\mu_s^{{X},N})\right]\leq \E\left[\big|e_s^i\big|^2\right].$$
By utilising Assumption \ref{assp1} and Young's inequality, we  obtain
\begin{align*}
	\E\Big[\big|e_t^i\big|^p\Big]\leq& C_{p,L}\int_{0}^{t}\E\Big[\big|e_s^i\big|^p\Big]ds\\
	&+C_p\int_{0}^{t}\E\left[\big(2+\tilde{Q}_L(\hat{X}_s^{i,N},\bar{X}_s^{i,N})^p+\tilde{Q}_K(\hat{X}_s^{i,N},\bar{X}_s^{i,N})^p\big)\big|\hat{X}_s^{i,N}-\bar{X}_s^{i,N}\big|^p\right]ds.
\end{align*}
 Using H\"{o}lder's inequality means
 \begin{align*}
 	\E&\Big[\big(2+\tilde{Q}_L(\hat{X}_s^{i,N},\bar{X}_s^{i,N})^p+\tilde{Q}_K(\hat{X}_s^{i,N},\bar{X}_s^{i,N})^p\big)\big|\hat{X}_s^{i,N}-\bar{X}_s^{i,N}\big|^p\Big]\\
 	&\leq\Big(\E\Big[\big(2+\tilde{Q}_L(\hat{X}_s^{i,N},\bar{X}_s^{i,N})^p+\tilde{Q}_K(\hat{X}_s^{i,N},\bar{X}_s^{i,N})^p\big)^2\Big]\E\Big[\big|\hat{X}_s^{i,N}-\bar{X}_s^{i,N}\big|^{2p}\Big]\Big)^{1/2}\\
 	&\leq\Big(\E\Big[4+\tilde{Q}_L(\hat{X}_s^{i,N},\bar{X}_s^{i,N})^{2p}+\tilde{Q}_K(\hat{X}_s^{i,N},\bar{X}_s^{i,N})^{2p}\Big]\E\Big[\big|\hat{X}_s^{i,N}-\bar{X}_s^{i,N}\big|^{2p}\Big]\Big)^{1/2}.
 \end{align*}
 Owing to the boundedness established in Theorem \ref{theorem1}, we arrive at $$\E\left[\tilde{Q}_K(\hat{X}_s^{i,N},\bar{X}_s^{i,N})^{2p}\right]\leq C,\qquad \E\left[\tilde{Q}_L(\hat{X}_s^{i,N},\bar{X}_s^{i,N})^{2p}\right]\leq C.$$
 From (\ref{3.3}) for any $s\in [0,T],$ applying H\"{o}lder's inequality yields
 \begin{align*}
 	\E\Big[\big|\hat{X}_s^{i,N}-\bar{X}_s^{i,N}\big|^{2p}\Big]&\leq3^{2p-1}\left(\E\Big[\big|b(\bar{X}_s^{i,N},\mu_s^{\bar{X},N})\big|^{4p}\Big]\E\Big[(s-\underline{s})^{4p}\Big]\right)^{1/2}\\
 	&\quad+3^{2p-1}\left(\E\Big[\big\|\sigma(\bar{X}_s^{i,N},\mu_s^{\bar{X},N})\big\|^{4p}\Big]\E\Big[\big|W_s^i-W_{\underline{s}}^i\big|^{4p}\Big]\right)^{1/2}\\
 	&\quad+3^{2p-1}\left(\E\Big[\big\|\sigma^0(\bar{X}_s^{i,N},\mu_s^{\bar{X},N})\big\|^{4p}\Big]\E\Big[\big|W_s^0-W_{\underline{s}}^0\big|^{4p}\Big]\right)^{1/2},
 \end{align*}
 where $\E\left[\big|b(\bar{X}_s^{i,N},\mu_s^{\bar{X},N})\big|^{4p}\right],$  $\E\Big[\big\|\sigma(\bar{X}_s^{i,N},\mu_s^{\bar{X},N})\big\|^{4p}\Big]$ and  $\E\Big[\big\|\sigma^0(\bar{X}_s^{i,N},\mu_s^{\bar{X},N})\big\|^{4p}\Big]$ are uniformly bounded on $[0,T]$ by Remark \ref{rmk1} and Theorem \ref{theorem1}. Additionally, we have $$\E\big[(s-\underline{s})^{4p}\big]\leq(\delta T)^{4p}\leq C\delta^{2p},$$ $$\E\Big[\big|W_s^i-W_{\underline{s}}^i\big|^{4p}\Big]=\E\left[\E\Big[\big|W_s^i-W_{\underline{s}}^i\big|^{4p}\Big]\Big|{\cal F}_{\underline{s}}\right]\leq c_p\E\big[(s-\underline{s})^{2p}\big]\leq C\delta^{2p},$$ $$\E\Big[\big|W_s^0-W_{\underline{s}}^0\big|^{4p}\Big]=\E\left[\E\Big[\big|W_s^0-W_{\underline{s}}^0\big|^{4p}\Big]\Big|{\cal F}_{\underline{s}}\right]\leq c_p\E\big[(s-\underline{s})^{2p}\big]\leq C\delta^{2p}.$$
Therefore, collecting all the estimates together, we get
 \begin{align*}
 	\max_{i\in\{1,\dots,N\}}\sup_{0\leq s\leq t}\E\Big[|e_s^i|^p\Big]\leq C\int_{0}^{t}\max_{i\in\{1,\dots,N\}}\sup_{0\leq u\leq s}\E\Big[|e_u^i|^p\Big]ds+C\delta^{p/2}.
 \end{align*}
The result can be obtained immediately by applying Gronwall's inequality.
 \eproof

\subsection{Boundedness and convergence in infinite horizon}\indent
\indent In this subsection, we  extend our analysis from finite time to infinite time by giving moment bounds and strong convergence order under some stronger assumptions. We modify Assumption \ref{assp1}(1) as  the following  assumption:
\begin{assp}\label{assp3}
	There exist constants $\lambda_1,\lambda_2>0$ such that
	\begin{align*}
		\big\langle& x-x{'},b(x,\mu)-b(x{'},\mu{'})\big\rangle+(p-1)\|\sigma(x,\mu)-\sigma(x{'},\mu{'}) \|^2\\
		&+(p-1)\|\sigma^0(x,\mu)-\sigma^0(x{'},\mu{'}) \|^2\le -\lambda_1| x - x{'} |^2+\lambda_2\mathcal{W}_2^2(\mu,\mu{'}),
	\end{align*} 
	for any $x,x{'}\in \mathbb R^d $ and $\ \mu,\mu{'}\in {\cal P}_2(\mathbb R^d)$.
\end{assp}

\begin{rmk}\label{rmk2}
	Based on Assumption \ref{assp3}, for any $\varepsilon>0$, we obtain 
	\begin{align}\label{R3.12}
		\big\langle& x,b(x,\mu)\big\rangle+\frac{p-1}{2}\|\sigma(x,\mu)\|^2+\frac{p-1}{2}\|\sigma^0(x,\mu)\|^2\nonumber\\
		\leq&\big\langle x,b(x,\mu)-b(0,\delta_0)\big\rangle+\big\langle x,b(0,\delta_0)\big\rangle+(p-1)\|\sigma(x,\mu)-\sigma(0,\delta_0)\|^2+(p-1)\|\sigma(0,\delta_0)\|^2\nonumber\\
		&+(p-1)\|\sigma^0(x,\mu)-\sigma^0(0,\delta_0)\|^2+(p-1)\|\sigma^0(0,\delta_0)\|^2\nonumber\\
		\leq&-\lambda_1|x|^2+\lambda_2\mathcal{W}_2^2(\mu,\delta_0)+|x||b(0,\delta_0)|+(p-1)\|\sigma(0,\delta_0)\|^2+(p-1)\|\sigma^0(0,\delta_0)\|^2\nonumber\\
		\leq&-\lambda_1|x|^2+\lambda_2\mathcal{W}_2^2(\mu,\delta_0)+\varepsilon|x|^2+\frac{1}{\varepsilon}|b(0,\delta_0)|^2+(p-1)\|\sigma(0,\delta_0)\|^2+(p-1)\|\sigma^0(0,\delta_0)\|^2\nonumber\\
		=&(\varepsilon-\lambda_1)|x|^2+\lambda_2\mathcal{W}_2^2(\mu,\delta_0)+\frac{1}{\varepsilon}|b(0,\delta_0)|^2+(p-1)\|\sigma(0,\delta_0)\|^2+(p-1)\|\sigma^0(0,\delta_0)\|^2.
	\end{align}
	Choose sufficiently small $\varepsilon$ to make $-\lambda_1+\varepsilon<0$ hold. Then from (\ref{R3.12}), we  see that there exist constants $\gamma_1>0,\eta>0$ such that
	\begin{align*}
		\big\langle x,b(x,\mu)\big\rangle+\frac{p-1}{2}\|\sigma(x,\mu)\|^2+\frac{p-1}{2}\|\sigma^0(x,\mu)\|^2\leq -\gamma_1|x|^2+\lambda_2{\cal W}_2^2(\mu,\delta_0)+\eta,
	\end{align*}
	where $\gamma_1=\lambda_1-\varepsilon$, $\eta=\frac{1}{\varepsilon}|b(0,\delta_0)|^2+(p-1)\|\sigma(0,\delta_0)\|^2+(p-1)\|\sigma^0(0,\delta_0)\|^2.$ 
\end{rmk}
\begin{assp}\label{assp4}(Adaptive timestep for infinite time interval)\\\\
(1) The adaptive timestep function $h:\RR^d \rightarrow \RR^+$ is continuous and bounded by $h_{\max}$ with $0<h_{\max}<\8$, and there exists a positive constant $C_h$ such that for all $x\in \RR^d$ and $\mu\in {\cal P}_2(\mathbb R^d),$ $h$ satisfies 
\begin{equation*}
	h(x,\mu)= \frac{C_h}{\big(1+|b(x,\mu)|\|\sigma(x,\mu)\|+|b(x,\mu)|\|\sigma^0(x,\mu)\|+|x|^q\big)^{2}},
\end{equation*}
where $q$ is the same as that in Assumption \ref{assp1}.\\
(2) The timestep function $h^\delta:\RR^d \rightarrow \RR^+,0<\delta \leq1$ satisfies
\begin{equation*}
		\delta \min\big(h_{\max},h(x,\mu)\big)\leq h^\delta(x,\mu)\leq \min \big(\delta h_{\max},h(x,\mu)\big).
\end{equation*}
\end{assp}

We can choose $h(x,\mu)$ appropriately so that it satisfies $|b(x,\mu)|h(x,\mu)\leq C(1+|x|+{\cal W}_2(\mu,\delta_0))$ and $\|\sigma(x,\mu)\|h(x,\mu)^{\frac{1}{2}}+\|\sigma^0(x,\mu)\|h(x,\mu)^{\frac{1}{2}}\leq C(1+|x|+{\cal W}_2(\mu,\delta_0))$.
The following result can be found in Proposition 2.5 in \cite{N24}, so we omit it here.

\begin{lemma}\label{lemma2}
	Let $X_t$ be a solution to McKean-Vlasov SDE (\ref{eq1.2}). If (\ref{eq1.2}) satisfies Assumption \ref{assp3}, then for all $r\in(0,\tilde{p})$ there exists a constant $C_p>0$ which is independent of $t$. For any $t\geq0$, we have $$E\left[\big|{X}_t\big|^r\right]\leq C_p.$$
\end{lemma}

\begin{theorem}(Moment bound in infinite time interval)\label{theorem3.11}
	If (\ref{eq1.2}) satisfies Assumption \ref{assp1}(2), Assumption \ref{assp3}, and the timestep function $h$ satisfies Assumption \ref{assp4}, then for all $2\leq p<\tilde{p},$ $\gamma_1>\lambda_2+1$, there exists a constant $C$ that is independent of $t,\delta$ and $N$, but is solely dependent on initial condition  $X_0,$ the exponent $p,$ and parameters $\eta,\gamma_1,\lambda_2 $ such that for all $t\geq0$,
	$$\E\left[\big|\hat{X}_t^{i,N}\big|^p\right]\leq C.$$
\end{theorem}\indent
\indent\proof
First, for some $C_p\geq0$, by using (\ref{3.3}) we get
\begin{align}\label{3.10}
	\E\Big[\big|\hat{X}_s^{i,N}\big|^p\Big]\leq&4^{p-1}\E\Big[\big|\bar{X}_s^{i,N}\big|^p+\big|b(\bar{X}_{s}^{i,N},\mu_{s}^{\bar{X},N})\big|^p\big|s-\underline{s}\big|^p+\big|\sigma(\bar{X}_{s}^{i,N},\mu_{s}^{\bar{X},N})(W_s^i-W_{\underline{s}}^i)\big|^{p}\nonumber\\
	&+\big|\sigma^0(\bar{X}_{s}^{i,N},\mu_{s}^{\bar{X},N})(W_s^0-W_{\underline{s}}^0)\big|^{p}\Big]\nonumber\\
	\leq&C_p\E\Big[1+\big|\bar{X}_s^{i,N}\big|^p+{\cal W}_2^p(\mu_s^{\bar{X},N},\delta_0)\big]+4^{p-1}\E\big[\E\big[\big|\sigma(\bar{X}_{s}^{i,N},\mu_{s}^{\bar{X},N})(W_s^i-W_{\underline{s}}^i)\big|^{p}\Big|{\cal F}_{\underline{s}}\Big]\Big]\nonumber\\
	&+4^{p-1}\E\Big[\E\Big[\big|\sigma^0(\bar{X}_{s}^{i,N},\mu_{s}^{\bar{X},N})(W_s^0-W_{\underline{s}}^0)\big|^{p}\Big|{\cal F}_{\underline{s}}\Big]\Big]\nonumber\\
	\leq&C_p\E\Big[1+\big|\bar{X}_s^{i,N}\big|^p\Big],
\end{align}
where $\E\left[\big|W_s-W_{\underline{s}}\big|^p\Big|{\cal F}_{\underline{s}}\right]\leq C_p(s-\underline{s})^{p/2}$ and $\E\Big[{\cal W}_2^p(\mu_{s}^{\bar{X},N},\delta_0)\Big]=\E\Big[\big|\bar{X}_s^{i,N}\big|^p\Big]$  are used.
Moreover, using Remark \ref{rmk2} and (\ref{3.3}) yields 
\begin{align}\label{3.11}
	\E\Big[&\big|\hat{X}_s^{i,N}\big|^{p-2}\big(\Big\langle\hat{X}_{s}^{i,N},b(\bar{X}_{s}^{i,N},\mu_{s}^{\bar{X},N})\Big\rangle+\frac{(p-1)}{2}\big\|\sigma(\bar{X}_{s}^{i,N},\mu_{s}^{\bar{X},N})\big\|^2+\frac{(p-1)}{2}\big\|\sigma^0(\bar{X}_{s}^{i,N},\mu_{s}^{\bar{X},N})\big\|^2\big)\Big]\nonumber\\
	\leq&C_p\E\Big[\big(\big|\hat{X}_s^{i,N}-\bar{X}_{s}^{i,N}\big|^{p-2}+\big|\bar{X}_{s}^{i,N}\big|^{p-2}\big)\big(\Big\langle\hat{X}_{s}^{i,N}-\bar{X}_{s}^{i,N},b(\bar{X}_{s}^{i,N},\mu_{s}^{\bar{X},N})\Big\rangle\nonumber\\
	&+\Big\langle\bar{X}_{s}^{i,N},b(\bar{X}_{s}^{i,N},\mu_{s}^{\bar{X},N})\Big\rangle+\frac{(p-1)}{2}\big\|\sigma(\bar{X}_{s}^{i,N},\mu_{s}^{\bar{X},N})\big\|^2+\frac{(p-1)}{2}\big\|\sigma^0(\bar{X}_{s}^{i,N},\mu_{s}^{\bar{X},N})\big\|^2\big)\Big]\nonumber\\
	\leq&C_p\E\Big[\big|\hat{X}_s^{i,N}-\bar{X}_{s}^{i,N}\big|^{p-2}\big(\Big\langle\bar{X}_{s}^{i,N},b(\bar{X}_{s}^{i,N},\mu_{s}^{\bar{X},N})\Big\rangle+\frac{(p-1)}{2}\big\|\sigma(\bar{X}_{s}^{i,N},\mu_{s}^{\bar{X},N})\big\|^2\nonumber\\
	&+\frac{(p-1)}{2}\big\|\sigma^0(\bar{X}_{s}^{i,N},\mu_{s}^{\bar{X},N})\big\|^2\big)\Big]+C_p\E\Big[\big|\bar{X}_{s}^{i,N}\big|^{p-2}\Big(\Big\langle\bar{X}_{s}^{i,N},b(\bar{X}_{s}^{i,N},\mu_{s}^{\bar{X},N})\Big\rangle\nonumber\\
	&+\frac{(p-1)}{2}\big\|\sigma(\bar{X}_{s}^{i,N},\mu_{s}^{\bar{X},N})\big\|^2+\frac{(p-1)}{2}\big\|\sigma^0(\bar{X}_{s}^{i,N},\mu_{s}^{\bar{X},N})\big\|^2\Big)\Big]\nonumber\\
	&+C_p\E\Big[\big(\big|\hat{X}_s^{i,N}-\bar{X}_{s}^{i,N}\big|^{p-2}+\big|\bar{X}_{s}^{i,N}\big|^{p-2}\big)\Big\langle\hat{X}_{s}^{i,N}-\bar{X}_{s}^{i,N},b(\bar{X}_{s}^{i,N},\mu_{s}^{\bar{X},N})\Big\rangle\Big]\nonumber\\
	\leq&C_p\E\Big[\big(1+\big|\bar{X}_{s}^{i,N}\big|^{p-2}\big)\big(-\gamma_1\big|\bar{X}_{s}^{i,N}\big|^2+\lambda_2{\cal W}_2^2(\mu_s^{\bar{X},N},\delta_0)+\eta\big)\Big]\nonumber\\
	&+C_p\E\Big[\big|\bar{X}_{s}^{i,N}\big|^{p-2}\big(-\gamma_1\big|\bar{X}_{s}^{i,N}\big|^2+\lambda_2{\cal W}_2^2(\mu_s^{\bar{X},N},\delta_0)+\eta\big)\Big]+C_p\E\Big[\big(1+\big|\bar{X}_{s}^{i,N}\big|^{p-2}\big)\big(1+\big|\bar{X}_{s}^{i,N}\big|^2\big)\Big]\nonumber\\
	\leq&C_p(\lambda_2-\gamma_1+1)\E\Big[\big|\bar{X}_{s}^{i,N}\big|^p\Big]+C_{p,\eta}.
\end{align}
Next, by applying It\^{o}'s formula to $e^{\gamma pt}\big|\hat{X}_t^{i,N}\big|^p$ for any $\gamma>0,$ we get 
\begin{align*}
	e^{\gamma pt}\big|\hat{X}_{t}^{i,N}\big|^p\leq&\big|\hat{X}_{0}^{i,N}\big|^p+p\int_0^t\gamma e^{\gamma ps}\big|\hat{X}_{s}^{i,N}\big|^pds+p\int_0^te^{\gamma ps}\big|\hat{X}_{s}^{i,N}\big|^{p-2}\Big(\left\langle\hat{X}_{s}^{i,N},b(\bar{X}_{s}^{i,N},\mu_{s}^{\bar{X},N})\right\rangle\\
	&+\frac{(p-1)}{2}\big\|\sigma(\bar{X}_{s}^{i,N},\mu_{s}^{\bar{X},N})\big\|^2+\frac{(p-1)}{2}\big\|\sigma^0(\bar{X}_{s}^{i,N},\mu_{s}^{\bar{X},N})\big\|^2\Big)ds\\
	&+p\int_0^te^{\gamma ps}\big|\hat{X}_{s}^{i,N}\big|^{p-2}\left\langle\hat{X}_{s}^{i,N},\sigma(\bar{X}_{s}^{i,N},\mu_{s}^{\bar{X},N}) dW_s^i\right\rangle\\
	&+p\int_0^te^{\gamma ps}\big|\hat{X}_{s}^{i,N}\big|^{p-2}\left\langle\hat{X}_{s}^{i,N},\sigma^0(\bar{X}_{s}^{i,N},\mu_{s}^{\bar{X},N}) dW_s^0\right\rangle.
\end{align*}
Then, it follows from (\ref{3.10}) and (\ref{3.11}) that
\begin{align*}
	\E\left[e^{\gamma pt}\big|\hat{X}_{t}^{i,N}\big|^p\right]\leq&\E\left[\big|\hat{X}_{0}^{i,N}\big|^p\right]+C_p\E\Big[\int_0^te^{\gamma ps}(\gamma+\lambda_2-\gamma_1+1)\big|\bar{X}_{s}^{i,N}\big|^pds\Big]+\E\Big[\int_0^te^{\gamma ps}C_{p,\gamma,\eta}ds\Big]\\
	\leq& \E\left[\big|\hat{X}_{0}^{i,N}\big|^p\right]+C_p(\gamma+\lambda_2-\gamma_1+1)\int_0^t\sup_{0\leq r\leq s}\E\Big[e^{\gamma pr}\big|\hat{X}_{r}^{i,N}\big|^p\Big]ds\\
	&+\E\Big[\int_0^te^{\gamma ps}C_{p,\gamma,\eta}ds\Big].
\end{align*}
We choose $\gamma=\gamma_1-\lambda_2-1>0$ gives
\begin{align*}
		\E\left[e^{\gamma pt}\big|\hat{X}_{t}^{i,N}\big|^p\right]\leq\E\left[\big|\hat{X}_{0}^{i,N}\big|^p\right]+e^{\gamma pt}C_{p,\gamma,\eta}+C_{p,\gamma,\eta}.
\end{align*}
Hence,
\begin{equation*}
	\E\left[\big|\hat{X}_{t}^{i,N}\big|^p\right]\leq C_{X_0,\gamma_1,\lambda_2,\eta,p}.
\end{equation*}
\eproof

Combining Theorem 4.1 in \cite{L23} and Proposition 3.2 in \cite{N24} gives the following results.
\begin{prop}\label{prop3}
	Let Assumption \ref{assp1}, Assumption \ref{assp3} be satisfied with $p\in[2,r)$, if $2\lambda_2<\lambda_1,$ and $\lambda_2<\gamma_1,$ then for some $\lambda\in(0,\gamma_1-\lambda_2)$ 
	\begin{equation*}
		\sup_{i\in\{1,\dots,N\}}\E\Big[\big|X_t^i-X_t^{i,N}\big|^2\Big]\le Ce^{-\lambda t}\left\{
		\begin{array}{l}
			N^{-1/2}+N^{-(r-2)/r},\quad\quad \quad\quad\quad if\  d<4\ and\ r\neq4,\\
			N^{-1/2}\log(1+N)+N^{-(r-2)/r},\: if\  d=4\ and\ r\neq4,\\
			N^{-2/d}+N^{-(r-2)/r},\quad \quad\quad\quad\quad  if\  d>4\ and\ r\neq d/(d-2),
		\end{array}\right. 
	\end{equation*}
	where the constant $C>0$ does not depend on $N$ and $T$.
\end{prop}

\begin{theorem}(Strong convergence order)
Let $p>0$ and $X_0\in L_0^p(\RR^d).$ If (\ref{eq1.2}) satisfies Assumption \ref{assp1}, Assumption \ref{assp3}, and the time-step function h satisfies Assumption \ref{assp4}, with $\lambda_1>\lambda_2+\frac{5}{2}$, then there exists a constant $C>0$ such that 
\begin{equation*}
	\max_{i\in \{1,\dots,N\}}\E\left[\left|\hat{X}_t^{i,N}-X_t^{i,N}\right|^p\right]\leq C\delta^{p/2}.
\end{equation*}
\end{theorem}\indent
\indent\proof
Define $e_t^i=\hat{X}_t^{i,N}-X_t^{i,N},$ then we have
\begin{align*}
	de_t^i=&\Big(b(\bar{X}_t^{i,N},\mu_t^{\bar{X},N})-b({X}_t^{i,N},\mu_t^{{X},N})\Big)dt+\Big(\sigma(\bar{X}_t^{i,N},\mu_t^{\bar{X},N})-\sigma({X}_t^{i,N},\mu_t^{{X},N})\Big)dW_t^i\\
	&+\Big(\sigma(\bar{X}_t^{i,N},\mu_t^{\bar{X},N})-\sigma({X}_t^{i,N},\mu_t^{{X},N})\Big)dW_t^0.
\end{align*}
For any $\lambda>0,$ by It\^{o}'s formula we obtain
\begin{align}\label{3.12}
	e^{{p\lambda t}}\big|e_t^i\big|^p\leq&\int_{0}^{t}p\lambda e^{p\lambda s}\big|e_s^i\big|^pds\nonumber\\
	&+\int_{0}^{t}pe^{p\lambda s}\big|e_s^i\big|^{p-2}\Big\langle e_s^i,b(\bar{X}_s^{i,N},\mu_s^{\bar{X},N})-b({X}_s^{i,N},\mu_s^{{X},N})\Big\rangle ds\nonumber\\
	&+\int_{0}^{t}\frac{p(p-1)}{2}e^{p\lambda s}\big|e_s^i\big|^{p-2}\big\|\sigma(\bar{X}_s^{i,N},\mu_s^{\bar{X},N})-\sigma({X}_s^{i,N},\mu_s^{{X},N})\big\|^2ds\nonumber\\
	&+\int_{0}^{t}\frac{p(p-1)}{2}e^{p\lambda s}\big|e_s^i\big|^{p-2}\big\|\sigma^0(\bar{X}_s^{i,N},\mu_s^{\bar{X},N})-\sigma^0({X}_s^{i,N},\mu_s^{{X},N})\big\|^2ds\nonumber\\
	&+\int_{0}^{t}pe^{p\lambda s}\big|e_s^i\big|^{p-2}\Big\langle e_s^i,\big(\sigma(\bar{X}_s^{i,N},\mu_s^{\bar{X},N})-\sigma({X}_s^{i,N},\mu_s^{{X},N})\big) dW_s^i\Big\rangle\nonumber\\
	&+\int_{0}^{t}pe^{p\lambda s}\big|e_s^i\big|^{p-2}\Big\langle e_s^i,\big(\sigma^0(\bar{X}_s^{i,N},\mu_s^{\bar{X},N})-\sigma^0({X}_s^{i,N},\mu_s^{{X},N})\big) dW_s^0\Big\rangle.
\end{align}
Note that
\begin{align*}
	&\Big\langle e_s^i,b(\bar{X}_s^{i,N},\mu_s^{\bar{X},N})-b({X}_s^{i,N},\mu_s^{{X},N})\Big\rangle=\Big\langle e_s^i,b(\bar{X}_s^{i,N},\mu_s^{\bar{X},N})-b(\hat{X}_s^{i,N},\mu_s^{\hat{X},N})\Big\rangle\\
	&+\Big\langle e_s^i,b(\hat{X}_s^{i,N},\mu_s^{\hat{X},N})-b({X}_s^{i,N},\mu_s^{{X},N})\Big\rangle.
\end{align*}
Using Assumption \ref{assp1} and Young's inequality, it can be deduced that
\begin{align}\label{3.13}
	\Big\langle e_s^i,&b(\bar{X}_s^{i,N},\mu_s^{\bar{X},N})-b(\hat{X}_s^{i,N},\mu_s^{\hat{X},N})\Big\rangle\nonumber\\
	\leq&\big|e_s^i\big|\Big(\tilde{Q}_L\big(\hat{X}_s^{i,N},\bar{X}_s^{i,N}\big)\big|\hat{X}_s^{i,N}-\bar{X}_s^{i,N}\big|+{\cal W}_2(\mu_s^{\hat{X},N},\mu_s^{\bar{X},N})\Big)\nonumber\\
	\leq&\frac{1}{2}\big|e_s^i\big|^2+\tilde{Q}_L\big(\hat{X}_s^{i,N},\bar{X}_s^{i,N}\big)^2\big|\hat{X}_s^{i,N}-\bar{X}_s^{i,N}\big|^2+{\cal W}_2^2(\mu_s^{\hat{X},N},\mu_s^{\bar{X},N}),
\end{align}
where $\tilde{Q}_L(\hat{X}_s^{i,N},\bar{X}_s^{i,N}):=L(1+\big|\hat{X}_s^{i,N}\big|^q+\big|\bar{X}_s^{i,N}\big|^q).$  Likewise, applying Remark \ref{rmk1} and Young's inequality on $\sigma$ means
\begin{align*}
	\big\|&\sigma(\bar{X}_s^{i,N},\mu_s^{\bar{X},N})-\sigma({X}_s^{i,N},\mu_s^{{X},N})\big\|^2\\
	&\leq2\big\|\sigma(\bar{X}_s^{i,N},\mu_s^{\bar{X},N})-\sigma(\hat{X}_s^{i,N},\mu_s^{\hat{X},N})\big\|^2 +2\big\|\sigma(\hat{X}_s^{i,N},\mu_s^{\hat{X},N})-\sigma({X}_s^{i,N},\mu_s^{{X},N})\big\|^2.
\end{align*}
We see 
\begin{align}\label{3.14}
	\big\|&\sigma(\bar{X}_s^{i,N},\mu_s^{\bar{X},N})-\sigma(\hat{X}_s^{i,N},\mu_s^{\hat{X},N})\big\|^2\nonumber\\
	&\leq\tilde{Q}_K\big(\hat{X}_s^{i,N},\bar{X}_s^{i,N}\big)|\hat{X}_s^{i,N}-\bar{X}_s^{i,N}|^2+{\cal W}_2^2(\mu_s^{\hat{X},N},\mu_s^{\bar{X},N}),
\end{align}
where $\tilde{Q}_K(\hat{X}_s^{i,N},\bar{X}_s^{i,N}):=K(1+\big|\hat{X}_s^{i,N}\big|^q+\big|\bar{X}_s^{i,N}\big|^q).$  Similarly, one obtains
\begin{align*}
	\big\|&\sigma^0(\bar{X}_s^{i,N},\mu_s^{\bar{X},N})-\sigma^0({X}_s^{i,N},\mu_s^{{X},N})\big\|^2\\
	&\leq2\big\|\sigma^0(\bar{X}_s^{i,N},\mu_s^{\bar{X},N})-\sigma^0(\hat{X}_s^{i,N},\mu_s^{\hat{X},N})\big\|^2\\
	&\quad+2\big\|\sigma^0(\hat{X}_s^{i,N},\mu_s^{\hat{X},N})-\sigma^0({X}_s^{i,N},\mu_s^{{X},N})\big\|^2,
\end{align*}
where
\begin{align}\label{3.15}
	\big\|&\sigma^0(\bar{X}_s^{i,N},\mu_s^{\bar{X},N})-\sigma^0(\hat{X}_s^{i,N},\mu_s^{\hat{X},N})\big\|^2\nonumber\\
	&\leq\tilde{Q}_K\big(\hat{X}_s^{i,N},\bar{X}_s^{i,N}\big)\big|\hat{X}_s^{i,N}-\bar{X}_s^{i,N}\big|^2+{\cal W}_2^2(\mu_s^{\hat{X},N},\mu_s^{\bar{X},N}).
\end{align}
By Assumption \ref{assp3}, we have
\begin{align}\label{3.16}
	\Big\langle& e_s^i,b(\hat{X}_s^{i,N},\mu_s^{\hat{X},N})-b({X}_s^{i,N},\mu_s^{{X},N})\Big\rangle+{(p-1)}\big\|\sigma(\hat{X}_s^{i,N},\mu_s^{\hat{X},N})-\sigma({X}_s^{i,N},\mu_s^{{X},N})\big\|^2\nonumber\\
	&\quad+{(p-1)}\big\|\sigma^0(\hat{X}_s^{i,N},\mu_s^{\hat{X},N})-\sigma^0({X}_s^{i,N},\mu_s^{{X},N})\big\|^2\nonumber\\
	&\leq-\lambda_1\big|e_s^i\big|^2+\lambda_2{\cal W}_2^2(\mu_s^{\hat{X},N},\mu_s^{X,N}).
\end{align}
Hence, substituting the derived estimates (\ref{3.13}), (\ref{3.14}), (\ref{3.15}), (\ref{3.16}) into (\ref{3.12}) yields
\begin{align*}
	e^{{p\lambda t}}\big|e_t^i\big|^p\leq&\int_{0}^{t}p(\lambda-\lambda_1+\frac{1}{2})e^{p\lambda s}\big|e_s^i\big|^pds+\int_{0}^{t}pe^{p\lambda s}\big|e_s^i\big|^{p-2}\lambda_2{\cal W}_2^2(\mu_s^{\hat{X},N},\mu_s^{{X},N})ds\\
	&+\int_{0}^{t}2p(p-1)e^{p\lambda s}\big|e_s^i\big|^{p-2}\Big(\tilde{Q}_K(\hat{X}_s^{i,N},\bar{X}_s^{i,N})\big|\hat{X}_s^{i,N}-\bar{X}_s^{i,N}\big|^2+{\cal W}_2^2(\mu_s^{\hat{X},N},\mu_s^{\bar{X},N})\Big)ds\\
	&+\int_{0}^{t}pe^{p\lambda s}\big|e_s^i\big|^{p-2}\Big(\tilde{Q}_L(\hat{X}_s^{i,N},\bar{X}_s^{i,N})^2\big|\hat{X}_s^{i,N}-\bar{X}_s^{i,N}\big|^2+{\cal W}_2^2(\mu_s^{\hat{X},N},\mu_s^{\bar{X},N})\Big)ds\\
	&+\int_{0}^{t}pe^{p\lambda s}\big|e_s^i\big|^{p-2}\Big\langle e_s^i,\big(\sigma(\bar{X}_s^{i,N},\mu_s^{\bar{X},N})-\sigma({X}_s^{i,N},\mu_s^{{X},N})\big)dW_s^i\Big\rangle\\
	&+\int_{0}^{t}pe^{p\lambda s}\big|e_s^i\big|^{p-2}\Big\langle e_s^i,\big(\sigma^0(\bar{X}_s^{i,N},\mu_s^{\bar{X},N})-\sigma^0({X}_s^{i,N},\mu_s^{{X},N})\big) dW_s^0\Big\rangle.
\end{align*}
Recall
$${\cal W}_2^2(\mu_s^{\hat{X},N},\mu_s^{\bar{X},N})ds\leq\frac{1}{N}\sum_{j=1}^{N}\big|\hat{X}_s^{j,N}-\bar{X}_s^{j,N}\big|^2,$$ and 
$${\cal W}_2^2(\mu_s^{\hat{X},N},\mu_s^{{X},N})\leq\frac{1}{N}\sum_{j=1}^{N}\big|e_s^j\big|^2.$$
Therefore, taking the expectation of both sides yields
\begin{align*}
	\E\Big[e^{{p\lambda t}}\big|e_t^i\big|^p\Big]\leq&\E\left[\int_{0}^{t}p\big(\lambda+\lambda_2-\lambda_1+\frac{1}{2}\big)e^{p\lambda s}\big|e_s^i\big|^pds\right]\\
	&+\E\Big[\int_{0}^{t}2p(p-1)e^{p\lambda s}\big|e_s^i\big|^{p-2}\Big(1+\tilde{Q}_K\big(\hat{X}_s^{i,N},\bar{X}_s^{i,N}\big)\Big)\big|\hat{X}_s^{i,N}-\bar{X}_s^{i,N}\big|^2ds\Big]\\
	&+\E\left[\int_{0}^{t}pe^{p\lambda s}\big|e_s^i\big|^{p-2}\Big(1+\tilde{Q}_L\big(\hat{X}_s^{i,N},\bar{X}_s^{i,N}\big)^2\Big)\big|\hat{X}_s^{i,N}-\bar{X}_s^{i,N}\big|^2ds\right],
\end{align*}
By applying Young's inequality, we can further analyze
\begin{align*}
	\E\Big[e^{{p\lambda t}}\big|e_t^i\big|^p\Big]\leq&\E\left[\int_{0}^{t}p\big(\lambda+\lambda_2-\lambda_1+\frac{5}{2}\big)e^{p\lambda s}\big|e_s^i\big|^pds\right]\\
	&+C_{p}\E\left[\int_{0}^{t}e^{p\lambda s}\left(2+\tilde{Q}_K\big(\hat{X}_s^{i,N},\bar{X}_s^{i,N}\big)^{p/2}+\tilde{Q}_L\big(\hat{X}_s^{i,N},\bar{X}_s^{i,N}\big)^p\right)\big|\hat{X}_s^{i,N}-\bar{X}_s^{i,N}\big|^pds\right].
\end{align*}
Since $C$ is independent of $t$, we chose $\lambda=\lambda_1-\lambda_2-\frac{5}{2}>0,$ so
\begin{align*}
	\E\Big[e^{{p\lambda t}}\big|e_t^i\big|^p\Big]\leq C_{p}\int_{0}^{t}e^{p\lambda s}\E\Big[\left(2+\tilde{Q}_K\big(\hat{X}_s^{i,N},\bar{X}_s^{i,N}\big)^{p/2}+\tilde{Q}_L\big(\hat{X}_s^{i,N},\bar{X}_s^{i,N}\big)^p\right)\big|\hat{X}_s^{i,N}-\bar{X}_s^{i,N}\big|^p\Big]ds.
\end{align*}
By H\"{o}lder's inequality, we derive
\begin{align*}
	\E\Big[&\left(2+\tilde{Q}_K\big(\hat{X}_s^{i,N},\bar{X}_s^{i,N}\big)^{p/2}+\tilde{Q}_L\big(\hat{X}_s^{i,N},\bar{X}_s^{i,N}\big)^p\right)\big|\hat{X}_s^{i,N}-\bar{X}_s^{i,N}\big|^p\Big]\\
	&\leq\left(\E\Big[\left(2+\tilde{Q}_K\big(\hat{X}_s^{i,N},\bar{X}_s^{i,N}\big)^{p/2}+\tilde{Q}_L\big(\hat{X}_s^{i,N},\bar{X}_s^{i,N}\big)^p\right)^2\Big]\E\Big[\big|\hat{X}_s^{i,N}-\bar{X}_s^{i,N}\big|^{2p}\Big]\right)^{1/2}.
\end{align*}
By Theorem \ref{theorem3.11}, there exists a constant $C>0$ such that$$\E\Big[\tilde{Q}_L\big(\hat{X}_s^{i,N},\bar{X}_s^{i,N}\big)^{2p}\Big]\leq C,\qquad\E\Big[\tilde{Q}_K\big(\hat{X}_s^{i,N},\bar{X}_s^{i,N}\big)^{p}\Big]\leq C.$$
And by the same argument as Theorem \ref{theorem3}, it is known that $$\E\Big[\big|\hat{X}_s^{i,N}-\bar{X}_s^{i,N}\big|^{2p}\Big]\leq C\delta^{p}.$$
Hence, we have
$$\E\Big[e^{{p\lambda t}}\big|e_t^i\big|^p\Big]\leq C_{L,p,h_{\max}}\int_{0}^{t}\delta^{p/2}e^{p\lambda s}ds.$$
Then,
$$\E\Big[\big|e_t^i\big|^p\Big]\leq C_{L,p,h_{\max},\lambda_1,\lambda_2}\delta^{p/2}.$$
\eproof

\section{Numerical examples}\label{section4}\indent
\indent The purpose of this section is to demonstrate the performance of the scheme proposed in this paper. To approximate the law ${\cal L}_{X_{t_n}}^1$ at each timestep $t_n$ on the interval $[0,T]$, we use a standard particle method with $N$ particles for each realization of $W^0$. For our experiment, we use $N=10^4$ to approximate the conditional law. Since the exact solution of the example is not known, to illustrate the strong convergence in $h$, we therefore compute the root-mean-square error (RMSE) by comparing the numerical solution at level $l$ with the solution at level $l-1$, at the final time $T$. This is denoted by RMSE, $$RMSE:=\sqrt{\frac{1}{N}\sum_{i=1}^{N}\Big(\hat{X}_T^{i,N,M_l}-\hat{X}_T^{i,N,M_{l-1}}\Big)^2},$$ where $\hat{X}_T^{i,N,M_l}$ denotes the numerical solution of $X$ at time $T$ (for our numerical experiment set $T=1$). $M_l=\lceil2^lT\rceil$ and the two-particle systems at each level are generated by the same Brownian motions. In adaptive timestep function $h^\delta,$ we use $\delta=M_l^{-1}$ for some fixed $M>1.$
\begin{expl}
	Consider the one-dimensional McKean-Vlasov SDE
\begin{equation} \label{4.1}
dX_t=\left(X_t-8X_t^3+\frac{1}{2}\E^1[X_t]\right)dt+\frac{1}{2}\left(X_t^2+\E^1[X_t]\right)dW_t+\frac{1}{2}\left(X_t^2+\E^1[X_t]\right)dW_t^0,\quad t\in[0,1],
\end{equation}
with $X_0=\sin (W(0))$. So 
$$b(x,\mu)=x-8x^3+\frac{1}{2}\int_{\mathbb{R}}x\mu(dx),$$  $$\sigma(x,\mu)=\frac{1}{2}\left(x^2+\int_{\mathbb{R}}x\mu(dx)\right),\sigma^0(x,\mu)= \frac{1}{2}\left(x^2+\int_{\mathbb{R}}x\mu(dx)\right).$$
 The coefficients can be readily verified to satisfy all assumptions listed above, thus it has a unique strong solution.
Consider the particle system of (4.1) as 
\begin{align}
	dX_t^{i,N}=&\left(X_t^{i,N}-8(X_t^{i,N})^3+\frac{1}{2N}\sum_{j=1}^NX_t^{j,N}\right)dt+\frac{1}{2}\left((X_t^{i,N})^2+\frac{1}{N}\sum_{j=1}^NX_t^{j,N}\right)dW_t^i\nonumber\\
	&+\frac{1}{2}\left((X_t^{i,N})^2+\frac{1}{N}\sum_{j=1}^NX_t^{j,N}\right)dW_t^0.
\end{align}
For the adaptive step size, choose $$h(x,\mu)=(\frac{1}{1+8|x|^5+\frac{1}{2}\int_{\mathbb{R}}|x|^2\mu(dx)})^2,\qquad h(x,\mu)^\delta=\delta h(x,\mu).$$ Set $C_3=1$, then we have
 $$|x-8x^3+\frac{1}{2}\int_{\mathbb{R}}^{}x\mu(dx)||x^2+\int_{\mathbb{R}}^{}x\mu(dx)|h(x,\mu)^{\frac{1}{2}}\leq 1.$$
Hence, $h(x,\mu)\leq(\frac{1}{|x-8x^3+\frac{1}{2}\int_{\mathbb{R}}x\mu(dx)||x^2+\int_{\mathbb{R}}x\mu(dx)|})^2$.  It can  be easily verified that Assumption \ref{assp2} holds by choosing the coefficients appropriately. So the adaptive EM scheme can work on this example.
 \end{expl}
 \begin{expl}
 	Consider the one-dimensional  McKean-Vlasov SDE
 \begin{equation}
 	dX_t=\left(-2X_t-3X_t^2|X_t|-2\E^1[X_t]\right)dt+\frac{1}{4}\left(1+|X_t|^{1.5}+\E^1[X_t]\right)dW_t+\frac{1}{4}\left(|X_t|^{1.5}+\E^1[X_t]\right)dW_t^0,
 \end{equation}
with $X_0=\sin (W(0))$. So $$b(x,\mu)=-2x-3x^2|x|-2\int_{\mathbb{R}}x\mu(dx),$$   $$\sigma(x,\mu)=\frac{1}{4}\left(1+|x|^{1.5}+\int_{\mathbb{R}}x\mu(dx)\right),\sigma^0(x,\mu)=\frac{1}{4}\left(1+|x|^{1.5}+\int_{\mathbb{R}}x\mu(dx)\right),$$
for all $x\in \mathbb{R}$ and $\mu\in\mathcal{P}(\mathbb{R})$. It is easy to show that these coefficients satisfy Assumption \ref{assp3}.  Set $h_{\max}=1,$ then we  choose 
$$h(x,\mu)=\min(1,(3|x|^3+2\int_{\mathbb{R}}|x|^2\mu(dx))^{-1}),\qquad h^\delta(x,\mu)=\delta\min(1,h(x,\mu)).$$ 
These choices show that the adaptive EM scheme can work on this example.
\end{expl}

\section*{Funding}\indent

This work is supported by the National Natural Science Foundation of China (62373383, 62076106), the Fundamental Research Funds for the Central Universities, South-Central MinZu University (CZQ25020, CZZ25007), and the Fund for Academic Innovation Teams of South-Central Minzu University (XTZ24004).


\begin{thebibliography}{99} 
\bibitem{BJ20}
J.~Bao, C.~Reisinger, P.~Ren, W.~Stockinger, {\em Milstein schemes and antithetic multilevel Monte Carlo sampling for delay McKean–Vlasov equations and interacting particle systems}, IMA J. Numer. Anal., 44(4): 2437-2479,
\bibitem{B21}
J.~Bao, C.~Reisinger, P.~Ren, W.~Stockinger, {\em First-order convergence of Milstein schemes for McKean-Vlasov equations and interacting particle systems}, Proc. R. Soc. A., 477 (2245): 20200258(2021).
\bibitem{B18}
M.~Bauer, T.~Meyer-Brandis and F.~Proske, {\em Strong solutions of mean-field stochastic differential equations with irregular drift}, Electron. J. Probab., 23:1-35(2018).
\bibitem{B18}
D. Belomestny and J. Schoenmakers, {\em Projected particle methods for solving McKean–Vlasov stochastic differential equations}, SIAM J. Numer. Anal., 56(6):3169–3195(2018).
\bibitem{B97}
M.~Bossy and D.~Talay, {\em A stochastic particle method for the McKean–Vlasov and the Burgers equation}, Math. Comp., 66: 157–192(1997).
\bibitem{B22}
U.~Botija-Munoz, C.~Yuan, {\em Explicit Numerical Approximations for SDDEs in Finite and Infinite Horizons using the Adaptive EM Method: Strong Convergence and Almost Sure Exponential Stability}, Appl. Math. Comput., 478: 128853(2024).
\bibitem{C18}
R.~Carmona, and F.~Delarue, {\em Probabilistic Theory of Mean Field Games with Applications I:Mean Field FBSDEs,Control, and Games}, Probab. Theory Stoch. Model. 83(2018).
\bibitem{CR18}
R.~Carmona, and F.~Delarue, {\em Probabilistic Theory of Mean Field Games with Applications II: Mean Field Games with Common Noise and Master Equations}, Probab. Theory Stoch. Model. 84(2018).
\bibitem{R17}
G.~dos Reis, W.~Salkeld, J.~Tugaut, {\em Freidlin–Wentzell LDPs in path space for McKean-Vlasov equations and the functional iterated logarithm law}, Ann. Appl. Probab., 29: 1487-1540(2017).
\bibitem{D22}
G.~dos Reis, s.~Engelhardt, G.~Smith, {\em Simulation of McKean–Vlasov SDEs with super-linear growth}, IMA J. Numer. Anal., 42: 874-922(2022).
\bibitem{E23}
X.~Erny, E.~Löcherbach,D.~Loukianova, {\em Strong error bounds for the convergence to its mean field limit for systems of interacting neurons in a diffusive scaling}, Ann. Appl.Probab., 33(5): 3563-3586(2023).
\bibitem{W20}
W.~Fang, M.~Giles, {\em Adaptive Euler–Maruyama method for SDEs with nonglobally Lipschitz drift}, Ann. Appl.Probab., 30: 526–560(2020).
\bibitem{G24}
S.~Gao, Q.~Guo, J.~Hu, et al, {\em Convergence rate in Lp sense of tamed EM scheme for highly nonlinear neutral multiple-delay stochastic McKean–Vlasov equations}, J. Comput. Appl. Math., 441: 115682(2024).
\bibitem{HW21}
W.~Hammersley, D.~\v{S}i\v{s}ka, L.~Szpruch, {\em Weak existence and uniqueness for McKean–Vlasov SDEs with common noise}, Ann. Probab., 49: 527-555(2021).

\bibitem{H21}
W.~Hammersley, D.~\v{S}i\v{s}ka, L.~Szpruch, {\em Mckean-Vlasov SDE under measure dependent Lyapunov conditions}, Ann. Inst. Henri Poincaré Probab. Stat., 57: 1032-1057(2021).
\bibitem{Hu21}
X.~Huang, P.~Ren, F.~Wang, {\em Distribution dependent stochastic differential equations}, Front. Math. China, 16(2):257–301(2021).
\bibitem{H15}
M.~Hutzenthaler and A.~Jentzen, {\em Numerical approximations of stochastic differential equations with non-globally Lipschitz continuous coefficients}, Mem. Amer. Math. Soc., 236:1112(2015).
\bibitem{H12}
M.~Hutzenthaler, A.~Jentzen, and P.E.~Kloeden, {\em Strong convergence of an explicit numerical method for SDEs with nonglobally Lipschitz continuous coefficients}, Ann. Appl. Probab., 22: 1611–1641(2012).
\bibitem{Ka22}
A.~Kalinin, T.~Meyer-Brandis, F.~Proske, {\em Stability, uniqueness and existence of solutions to McKean-Vlasov SDEs in arbitrary moments}, J. Theoret. Probab., 37(4): 2941-2989(2024)
\bibitem{K92}
P.E.~Kloeden and E.~Platen, {\em Numerical Solution of Stochastic Differential Equations}, Springer Berlin Heidelberg(1992).
\bibitem{K22}
C.~Kumar, Neelima, C.~Reisinger, W.~Stockinger, {\em Well-posedness and tamed schemes for McKean-Vlasov equations with common noise}, Ann. Appl. Probab., 32: 3283-3330(2022).
\bibitem{K19}
C. Kumar and S. Sabanis, {\em  On Milstein approximations with varying coefficients: the case of super-linear diffusion coefficients}, BIT, 59(4): 929–968(2019).
\bibitem{L18}
D.~Lacker, {\em On a strong form of propagation of chaos for McKean–Vlasov equations}, Electron. Commun. Probab., 23(2018).
\bibitem{L21}
S.~Ledger and A.~S$\phi$jmark, {\em At the mercy of the common noise: blow-ups in a conditional McKean–Vlasov problem}, Electron. J. Probab., 26: 1-39(2021).
\bibitem{L22}
Y.~Li, X.~Mao, Q.~Song, F.~Wu, G.~Yin, {\em Strong convergence of Euler–Maruyama schemes for McKean-Vlasov stochastic differential equations under local Lipschitz conditions of state variables}, IMA J. Numer. Anal., 43(2): 1001-1035(2023).
\bibitem{L23}
Z.~Liu, S.~Gao, C.~Yuan, et al. {\em Stability of the numerical scheme for stochastic McKean-Vlasov equations}, arXiv preprint arXiv:2312.12699, (2023).
\bibitem{M66}
H.P.~McKean, {\em A class of Markov processes associated with nonlinear parabolic equations}, Proc. Natl. Acad. Sci. USA,  56 (6): 1907-1911(1966). 
\bibitem{Me20}
S.~Mehri, M.~Scheutzow, W.~Stannat, et al, {\em Propagation of chaos for stochastic spatially structured neuronal networks with delay driven by jump diffusions}, Ann. Appl. Probab., 30(1): 175-207(2020).
\bibitem{M20}
Y.S.~Mishura, and A.Y.~Veretennikov, {\em  Existence and uniqueness theorems for solutions of McKean–Vlasov stochastic equations}, Theory Probab. Math. Statist., 103: 59-101(2020). 
\bibitem{B20}
Neelima, S.~Biswas, C.~Kumar, G.~Reis and C.~Reisinger, {\em Well-posedness and tamed Euler schemes for McKean-Vlasov equations driven by L\'evy noise}, arXiv preprint \href{https://doi.org/10.48550/arXiv.2010.08585}{arXiv:2010.08585}(2020).
\bibitem{P16}
H.~Pham, {\em Linear quadratic optimal control of conditional McKean–Vlasov equation with random coefficients and applications}, Probab. Uncertain. Quant. Risk., 1, 1-26 (2016).
\bibitem{C22}
C.~Reisinger, W.~Stockinger, {\em An adaptive Euler–Maruyama scheme for McKean–Vlasov SDEs with super-linear growth and application to the mean-field FitzHugh–Nagumo model}, J. Comput. Appl. Math., 440: 113725(2022).
\bibitem{S13}
S.~Sabanis, {\em A note on tamed Euler approximations}, Electron. Commun. Probab., 18:1–10(2013).

\bibitem{S91}
A.S.~Sznitman, {\em Topics in propagation of chaos}, Ecole d’été de probabilités de Saint-Flour XIX — 1989. Lecture Notes in Math., 1464: 165--251(1991).
\bibitem{N24}
N.~K.~Tran, T.~T.~Kieu,D.~T.~Luong, et al. {\em On the infinite time horizon approximation for L\'evy-driven McKean-Vlasov SDEs with non-globally Lipschitz continuous and super-linearly growth drift and diffusion coefficients}, J. Math. Anal. Appl., 543(2): 128982(2025).

\bibitem{U18}
E.~Ullner, A.~Polito,and A.~Torcini, {\em Ubiquity of collective irregular dynamics in balanced networks of spiking neurons}, Chaos, 28(08)(2018).
















\end{thebibliography}
\end{document}